\documentclass[11pt,a4paper,reqno]{amsart}
\usepackage{bm}

\usepackage{setspace}
\usepackage{indentfirst}

\usepackage{tabls}
\usepackage{url}
\usepackage{color}
\definecolor{refkey}{rgb}{1,0,0}
\definecolor{labelkey}{rgb}{0,0,1}
\definecolor{labelkey}{rgb}{1,1,1}
\usepackage[linktocpage = true]{hyperref}
\hypersetup{
 colorlinks = true,
 citecolor = blue,
}

\usepackage{times}
\usepackage{amsthm}
\usepackage{amsmath}
\usepackage{amsfonts}
\usepackage{amssymb}
\usepackage{amsrefs}
\usepackage{mathrsfs}
\usepackage{xypic}
\usepackage{enumerate}
\usepackage{graphicx}
\usepackage{subfigure}
\usepackage[toc,page]{appendix}
\usepackage[all,cmtip]{xy}

\allowdisplaybreaks

\hypersetup{pdfpagemode=UseOutlines,colorlinks=false,pdfpagelayout=SinglePage,pdfstartview=FitH,bookmarksopen=true}

\newcommand{\abs}[1]{\lvert#1\rvert}

\newcommand{\prA}{\mathrm{Pr}_A}

\DeclareMathOperator{\CE}{CE}
\DeclareMathOperator{\Ha}{H}
\newtheorem{Thm}[equation]{Theorem}
\newtheorem*{main}{Theorem}
\newtheorem{Pro}[equation]{Proposition}
\newtheorem{Lem}[equation]{Lemma}

\newtheorem{thm-def}[equation]{Theorem-Definition}
\newtheorem{prop-def}[equation]{Proposition-Definition}
\newtheorem{def-prop}[equation]{Definition-Proposition}

\theoremstyle{definition}
\newtheorem{Def}[equation]{Definition}
\newtheorem{Ex}[equation]{Example}

\theoremstyle{remark}
\newtheorem{Rm}[equation]{Remark}
\numberwithin{equation}{section}
\begin{document}
\def\bp{\begin{proof}}
\def\ep{\end{proof}}
\def\be{\begin{equation}}
\def\ee{\end{equation}}
\def\bd{\begin{displaymath}}
\def\ed{\end{displaymath}}
\def\C{\mathcal{C}}
\def\D{\mathcal{D}}
\def\E{\mathcal{E}}
\def\F{\mathcal{F}}
\def\H{\textbf{H}}
\def\k{\mathbb{K}}
\def\G{\mathcal{G}}
\def\g{\mathfrak{g}}
\def\M{\mathcal{M}}
\def\N{\mathbb{N}}
\def\m{\mathfrak{m}}
\def\t{\mathfrak{t}}
\def\O{\mathcal{O}}
\def\P{\mathcal{P}}
\def\U{\mathcal{U}}
\def\V{\mathscr{V}}
\def\L{\mathcal{L}}
\def\X{\mathfrak{X}}
\def\Z{\mathbb{Z}}
\def\spec{\text{spec}}
\def\Im{\text{Im}}
\def\End{\operatorname{End}}
\def\Defo{\operatorname{Defo}}
\def\Der{\operatorname{Der}}
\def\Hom{\operatorname{Hom}}
\def\Jet{\operatorname{Jet}}
\def\Map{\operatorname{Map}}
\def\Mod{\operatorname{Mod}}
\def\sgn{\operatorname{sgn}}
\def\sh{\operatorname{sh}}

\title{Atiyah Classes of Strongly Homotopy Lie Pairs}
\thanks{Research partially supported by NSFC 11471179.}

\author{Zhuo Chen}%{$^\diamond$}\\
\address{Department of Mathematics, Tsinghua University}%, Beijing 100084, China}%\\\vspace{1mm} \\
\email{\href{mailto:chenzhuo@tsinghua.edu.cn}{chenzhuo@tsinghua.edu.cn}}

\author{Honglei Lang}%{$^\dag$}\\
\address{Department of Applied Mathematics, China Agricultural University}% 53072, Germany }  %\\\vspace{1mm}
\email{\href{mailto:hllang@cau.edu.cn}{hllang@cau.edu.cn}}

\author{Maosong Xiang}%$^*$\\
\address{Beijing International Center for Mathematical Research, Peking University}%, Beijing 100871, China} }
\email{\href{mailto:msxiang@pku.edu.cn}{msxiang@pku.edu.cn}}

%\footnotetext{{\it{Keywords}}:~Homotopical algebra, $L_\infty$-algebra, Atiyah class.}
%\footnotetext{{\it{MSC}}:~Primary  16E45, 18G55. Secondary   58C50.}
%\footnotetext{{\it{Emails}}:~\emph{$^\diamond$zchen@math.tsinghua.edu.cn,~~$^\dag$hllang@mpim-bonn.mpg.de,~~$^*$msxiang@pku.edu.cn}(corresponding author)}

%\maketitle
%\makeatletter

%\newif\if@borderstar

%\def\bordermatrix{\@ifnextchar*{%

%\@borderstartrue\@bordermatrix@i}{\@borderstarfalse\@bordermatrix@i*}%

%}

%\def\@bordermatrix@i*{\@ifnextchar{\@bordermatrix@ii}{\@bordermatrix@ii[()]}}

%\def\@bordermatrix@ii[#1]#2{%

%\begingroup

%\m@th\@tempdima8.75\p@\setbox\z@\vbox{%

%\def\cr{\crcr\noalign{\kern 2\p@\global\let\cr\endline }}%

%\ialign {$##$\hfil\kern 2\p@\kern\@tempdima & \thinspace %

%\hfil $##$\hfil && \quad\hfil $##$\hfil\crcr\omit\strut %

%\hfil\crcr\noalign{\kern -\baselineskip}#2\crcr\omit %

%\strut\cr}}%

%\setbox\tw@\vbox{\unvcopy\z@\global\setbox\@ne\lastbox}%

%\setbox\tw@\hbox{\unhbox\@ne\unskip\global\setbox\@ne\lastbox}%

%\setbox\tw@\hbox{%

%$\kern\wd\@ne\kern -\@tempdima\left\@firstoftwo#1%

%\if@borderstar\kern2pt\else\kern -\wd\@ne\fi%

%\global\setbox\@ne\vbox{\box\@ne\if@borderstar\else\kern 2\p@\fi}%

%\vcenter{\if@borderstar\else\kern -\ht\@ne\fi%

%\unvbox\z@\kern-\if@borderstar2\fi\baselineskip}%

%\if@borderstar\kern-2\@tempdima\kern2\p@\else\,\fi\right\@secondoftwo#1 $%

%}\null \;\vbox{\kern\ht\@ne\box\tw@}%

%\endgroup

%}

\makeatother

\begin{abstract}
  The subject of this paper is strongly homotopy (SH) Lie algebras, also known as $L_\infty$-algebras. We extract an intrinsic character, the Atiyah class, which measures the nontriviality of an SH Lie algebra $A$ when it is extended to $L$. In fact, given such an SH Lie pair $(L,A)$, and any $A$-module $E$, there associates a canonical cohomology class, the Atiyah class $[\alpha^E]$, which generalizes earlier known Atiyah classes out of Lie algebra pairs. We show that the Atiyah class $[\alpha^{L/A}]$ induces a graded Lie algebra structure on $\Ha^\bullet_{\mathrm{CE}}(A,L/A[-2])$, and the Atiyah class $[\alpha^E]$ of any $A$-module $E$ induces a Lie algebra module structure on $\Ha^\bullet_{\mathrm{CE}}(A,E)$. Moreover, Atiyah classes are invariant under gauge equivalent $A$-compatible infinitesimal deformations of $L$. 
\end{abstract}

\maketitle
\tableofcontents

\section*{Introduction}
This work is motivated by two sources: Atiyah classes and strongly homotopy Lie algebras.
Originally, the Atiyah class~\cite{Atiyah} of a holomorphic vector bundle $U$ over a complex manifold constitutes the obstruction to the existence of a holomorphic connection on $U$.  Molino~\cites{Molino1,Molino2} defined the Atiyah-Molino class of a foliation of a manifold to capture the existence of a locally projectable connection.
The Atiyah class of a Lie algebra pair was studied by Wang~\cite{Wang},  Nguyen-van \cite{Nguyen} and Bordemann~\cite{Bordemann}, to characterize the existence of invariant connections on a homogeneous space. %invariant connections in a principal bundle over a homogenous space.
%This problem was first done by Wang in \cite{Wang}, and then formulated as cohomological obstruction now called the Atiyah class by Nguyen-van \cite{Nguyen}. Bordemann \cite{Bordemann} generalized this result to the existence of equivariant connections on a homogeneous space.
Atiyah classes have enjoyed renewed vigor due to Kontsevich's seminal work on deformation quantization \cites{Kon1, Kon2}. %Kontsevich revealed the existence of deep ties between the Todd genus of complex manifolds and the Duflo element in Lie theory.
They are also related to the Rozansky-Witten theory~\cites{RW, Kap}.
%Among many recent works related to Atiyah classes, the paper ~\cite{CV} by
Calaque and Van den Bergh~\cite{CV} considered the Atiyah class of a DG module over a DG-algebra.  They also inferred that, given a  Lie algebra pair $(\mathfrak{d},\mathfrak{g})$,
the Atiyah class of the quotient $\mathfrak{d}/\mathfrak{g}$ coincides with the class capturing the obstruction to the ``PBW problem'' studied earlier by Calaque--C\u{a}ld\u{a}raru--Tu~\cite{CCT} (see also~\cites{Calaque,Grinberg}).

The notion of strongly homotopy (SH) Lie algebras, also called $L_\infty[1]$-algebras (see Definition \ref{Def:L[1]-algebra}), was introduced by Lada and Stasheff ~\cites{LS,LM}. The investigation of SH Lie algebras from various perspectives started a while ago. Attention on this subject in the past ten years is largely due to its role in mathematical physics and supergeometry.  For example, Kontsevich and Soibelman ~\cite{KS} approached this notion via the language of formal geometry.  Meanwhile, Bashkirov and Voronov~\cite{BV} used the Batalin-Vilkovisky formalism  to treat an SH Lie algebra as a special pointed BV$_\infty$-manifold. The interested reader is referred to a recent talk by Stasheff \cite{Stasheff}  {in which many related topics are reviewed.}

Motivated by the various constructions of Atiyah classes, we will study SH Lie algebra pairs, show that analogous Atiyah classes exist, and how they play the role of sending homotopical objects to Lie objects.
%SH Lie algebras are also called $L_\infty[1]$-algebras(see Definition \ref{Def:L[1]-algebra}).

It is usually nontrivial to construct examples of $L_\infty[1]$-algebras. One ``trivial'' way is the semi-direct product of an $L_\infty[1]$-algebra $A$ with its module $B$ (see Proposition~\ref{Thm:module} (4)). %as follows:~given an $L_\infty[1]$-algebra $A$ and an $A$-module $B$, one can endow $L=A\oplus B$ with a new $L_\infty[1]$-algebra structure, so called semi-direct product of the $L_\infty[1]$-algebra $A$ with its module $B$ (see Theorem \ref{Thm:module} (4)).
We are particularly interested in grasping the information when a smaller $L_\infty[1]$-algebra, say $A$, ``non-trivially'' extends to a bigger one, say $L$. Hence we introduce the SH Lie pair $(L,A)$, where $L$ is an $L_\infty[1]$-algebra and $A\subset L$ is a sub-algebra.
%A first notable fact is that
%As the quotient space $L/A$ is canonically an $A$-module (see Lemma \ref{Lem:AmodonB}), the semi. However, only knowing the $A$-module structure on $L/A$ is far from enough to determine the $L_\infty[1]$-structure underlying $L$, %although the semi-direct product of $A$ and $L/A$ exists and $L\cong A\oplus L/A$ (only as vector spaces).
%The semi-direct product $A\oplus L/A$ is certainly not interesting because  it does not tell how complicated $A$ and $L/A$ could be entangled in $L$.
What we discovered, is the so-called Atiyah class $[\alpha^{L/A}]$ of $(L,A)$, which generalizes previous constructions of Lie algebra pairs. It measures how ``nontrivial" it is  when the sub-algebra $A$ is extended to $L$, while the $A$-module structure on $L/A$ is maintained. Moreover, it refines the homotopical data of $(L,A)$, to a canonical graded Lie algebra $\Ha^\bullet_{\mathrm{CE}}(A,L/A[-2])$ (see Theorem~\ref{Cor:LieandLieModuleStr}).  One can even involve an external object --- an $A$-module $E$, and use the Atiyah class $[\alpha^E]$ of $E$ to test the nontrivial information of $A$ being extended to $L$. %Again, one shall obtain an
Moreover, $[\alpha^E]$ gives rise to a Lie algebra module structure on $\Ha^\bullet_{\mathrm{CE}}(A,E)$, over the aforesaid Lie algebra object (see Theorem~\ref{Cor:LieandLieModuleStr}).

%We shall give various characterizations of Atiyah classes and their properties, which will be followed by a number of low dimensional examples.
The following is a summary of this paper:

%Section~\ref{Preliminaries}  %is a succinct account of standard conventions of
%$\Z$-graded linear algebra and SH Lie algebras. %We will give and prove several equivalent descriptions of $L_\infty[1]$-modules and show that the category of such modules is in fact abelian.

After a review of $\Z$-graded linear algebra and SH Lie algebras in Section~\ref{Preliminaries}, we focus on the construction of Atiyah classes in Section~\ref{Atiyah classes}. Given an SH Lie pair $(L,A)$ and an $A$-module $E$, one is able to extend the $A$-module structure on $E$ to an $L$-connection $\nabla$ on $E$. The curvature $R^\nabla$ measures the failure of $E$ being an $L$-module. From $R^\nabla$, we extract a particular element
$$
\alpha_{\nabla}^E  := (J \otimes 1)(R^{\nabla}) \in \O(A) \otimes A^\perp \otimes \End(E),
$$
which is a degree $2$ cocycle. Here $\O(A)$ is the graded algebra of formal power series on $A$. We call $\alpha^E_\nabla$ the Atiyah cocycle of the SH Lie pair $(L,A)$ with respect to the $A$-module $E$ and the $L$-connection $\nabla$ extending $(A,E)$. %We then show the main property of such Atiyah cocycles:

\begin{main}[A]
The cohomology class, called Atiyah class, $[\alpha_\nabla^E]\in  \Ha^2_{\mathrm{CE}}(A,A^\perp\otimes \End(E))$ is canonical, i.e., independent of the choice of $\nabla$. In particular, for the canonical $A$-module $L/A$, there associates a canonical Atiyah class $[\alpha^{L/A}]\in \Ha^2_{\mathrm{CE}}(A,\Hom(L/A\otimes L/A,L/A))$.
\end{main}

In Section~\ref{Atiyah as functors}, we introduce the Atiyah operator and functor, which manifest the nature of Atiyah classes from different perspectives. The Atiyah operator is an $\O(A)$-linear map which arises from the construction of Atiyah cocycles:
\begin{align*}
&\bm{\alpha}^E:\qquad  \O(A)\otimes E \longrightarrow \O(A)\otimes A^\perp\otimes E;&\\
\mbox{Or,~}\quad&\bm{\alpha}^E:\qquad  (\O(A)\otimes L/A)~\times~(\O(A)\otimes E) \longrightarrow \O(A)\otimes E.&
%\overrightarrow{\alpha_\nabla^E}:\qquad & \O(A)\otimes B\longrightarrow \Hom_{\O(A)}(\O(A)\otimes E,\O(A)\otimes E)\cong \O(A)\otimes \End(E),
\end{align*}

We prove% The main result is the following

\begin{main}[B]
 The graded vector space $\Ha^\bullet_{\mathrm{CE}}(A,(L/A)[-2])$ with the binary operation induced by the Atiyah operator $\bm{\alpha}^{L/A}$ is a Lie algebra. Furthermore, if $E$ is an $A$-module, then $\Ha^\bullet_{\mathrm{CE}}(A,E)$  is a Lie algebra module over $\Ha_{\mathrm{CE}}^\bullet(A,(L/A)[-2])$, with the action induced by the Atiyah operator $\bm{\alpha}^E$.
\end{main}

This certainly generalizes previous results in \cites{Kap,CSX,Bottacin}. An alternative point of view is that the process of taking Atiyah classes defines a functor, called Atiyah functor, from the category of $A$-modules to the category of  $\Ha^\bullet_{\mathrm{CE}}(A,(L/A)[-2])$-modules.

In Section~\ref{Invariance section}, we study a special kind of deformations of the given SH Lie algebra pair $(L,A)$, namely  $A$-compatible infinitesimal deformations of $L$. Roughly speaking, they are $L_\infty[1]$-algebra structures on
$$
L[\hbar]=L\oplus \hbar L,
$$
where $\hbar$ is a formal parameter with $\hbar^2=0$, such that the subspace $A[\hbar]$ is trivially extended from $A$, and the $A[\hbar]$-module structure on $L[\hbar]/A[\hbar]=(L/A)[\hbar]$ is trivially extended from the $A$-module $L/A$.
A small perturbation of $L[\hbar]$ is an isomorphism $\sigma:\O(L)[\hbar]\rightarrow \O(L)[\hbar]$ of graded algebras. % We wish to see how Atiyah classes are affected when $L[\hbar]$ is perturbed.

We prove the following invariance property of Atiyah classes under gauge equivalences (See Definition~\ref{Def:gaugeequivalence}):

\begin{main}[C]
If two $A$-compatible infinitesimal deformations of $L$ are gauge equivalent, then the two associated Atiyah classes are the same.
\end{main}

It is our hope that these results may lead to new insights in homotopical algebras and DG-manifolds.
We would also like to point out other works that are related to the present paper:
Chen, Sti\'{e}non and Xu \cite{CSX} proposed a notion of the Atiyah class of a Lie algebroid pair $(L,A)$, which encompasses both the original Atiyah class of holomorphic vector bundles and the Atiyah-Molino class of a foliation as special cases. Shortly after that, an $L_\infty$-algebra structure on the space $\Gamma(\wedge^\bullet A^\vee \otimes L/A)$ was constructed in \cites{LSX1,LSX2}, where the Atiyah class determines the $2$-bracket $l_2$.
A similar theory for Lie groupoid pairs is available in \cite{LV}.
We also mention that Shoikhet~\cite{Shoikhet} studied the Atiyah class of a DG-manifold;
Costello~\cite{Costello} defined the Atiyah class of a DG-vector bundle in his geometric approach to Witten genus;
Mehta, Sti{\'e}non and Xu~\cite{MSX} studied the Atiyah class of a DG-Lie algebroid with respect to a DG-vector bundle.
%For more on this topic, we refer readers to the work~\cite{LV} of Laurent-Gengoux and Voglaire.

\section{Preliminaries}\label{Preliminaries}
\subsection{Graded linear algebra}
Throughout this paper, we fix a base field $\k$ of characteristic zero. A $\Z$-graded vector space is a  $\k$-vector space $V = \oplus_{n \in \Z} V^n$, where each $V^n = \{v \in V \mid \abs{v} = n\}$ is an ordinary $\k$-vector space consisting of elements of homogeneous degree $n$. Henceforth, we will simply call $V$ a graded vector space. And $\k$ is considered as concentrated in degree $0$.

A degree $r$ morphism from a graded vector space $V$ to a graded vector space $W$ is a linear map from $V$ to $W$ that sends $V^n$ to $W^{n+r}$, where $r$ could be any integer. The set $\Hom(V,W)$ consisting of such homogeneous morphisms is also a graded vector space. Thus the category of graded vector spaces over $\k$, denoted by $\text{GVS}_\k$, is a $\k$-linear category.

The dual of $V$, denoted by $V^\vee$, is the graded vector space whose degree $n$ part is the ordinary dual $(V^{-n})^\ast$ of $V^{-n}$. If $V$ is of finite dimension, then the dual of $V^\vee$ is isomorphic to $V$. In this paper, we will always assume that $V$ is finite dimensional if $V^\vee$ is involved.

For $k\in\Z$, we denote by $V[k]$ the graded vector space with $k$-shifted gradings $(V[k])^n = V^{n+k}$. Hence $(V[k])^\vee=V^\vee[-k]$.

The category of graded vector spaces is monoidal. The tensor product of two objects $V$ and $W$ is the graded vector space whose degree $n$ part is
\bd
(V \otimes W)^n  = \oplus_{i+j=n}V^i \otimes W^j.
\ed

We have isomorphisms of graded vector spaces:
\begin{align*}
  V^\vee \otimes W &\cong W\otimes V^\vee \cong \Hom(V,W),\\
  \xi\otimes w &\mapsto (-1)^{\abs{\xi}\abs{w}}w\otimes \xi\mapsto \phi(-),
\end{align*}
where $\xi\in V^\vee,w\in W$ and $\phi$ is the map
$v\mapsto (\xi\otimes w)(v) =(-1)^{\abs{v}\abs{w}} \xi(v)w$.

For any homogeneous element $\phi \in \Hom(V,W)$, its dual $\phi^\vee\in \Hom(W^\vee,V^\vee)$, which is also homogeneous of degree $\abs{\phi}$, is defined in the standard manner:
$$
\langle \phi^\vee(\alpha),v\rangle=(-1)^{\abs{\phi}\abs{\alpha}}
\langle \alpha,\phi(v)\rangle,\qquad \alpha\in W^\vee,v\in V.
$$

The symmetric algebra of $V$ and its formal completion are,  respectively,
\begin{align*}
S^\bullet(V)&=  \oplus_{n \geq 0} S^n(V), &  \widehat{S}{^\bullet(V)} &=  \prod_{n \geq 0}S^n(V).
\end{align*}
Note that they might be infinite dimensional.
The product in $S^\bullet(V)$, as well as that in $ \widehat{S}{^{\bullet}(V)}$, is denoted by $\odot$.
The Koszul sign $\epsilon{(\sigma)}$ of a permutation $\sigma$ of homogeneous vectors $v_1,\cdots,v_n$ in $V$ is determined by the equality
\bd
v_1 \odot \cdots \odot v_n = \epsilon(\sigma)v_{\sigma(1)} \odot \cdots \odot v_{\sigma(n)}.
\ed

Given $v \in V$, there induces two natural contractions, one from left and one from right, denoted respectively $\iota_v$ and $\llcorner v$, on $V^\vee$:
\bd
\iota_v \xi = (-1)^{\abs{\xi}\abs{v}}\xi \llcorner v = (-1)^{\abs{\xi}\abs{v}}\xi(v),\;\;\forall \xi\in V^\vee.
\ed
The left contraction $\iota_v$ is extended to
 $\iota_v:~S^\bullet(V^\vee) \rightarrow S^{\bullet-1}(V^\vee)$ by the Leibniz rule
\bd
\iota_v(\xi\odot\eta) = \iota_v(\xi)\odot\eta + (-1)^{\abs{v}\abs{\xi}}\xi\odot\iota_v(\eta),\;\;\forall \xi,\eta \in S^\bullet(V^\vee).
\ed
The extension of the right contraction is similar:
\bd
 (\xi\odot\eta)\llcorner v =  (-1)^{\abs{v}\abs{\eta}}(\xi\llcorner v)\odot\eta + \xi\odot (\eta\llcorner v),\;\;\forall \xi,\eta \in S^\bullet(V^\vee).
\ed

We define a duality pairing
\bd
S^\bullet(V)\times S^\bullet(V^\vee)\rightarrow \k
\ed
by
$$
\langle v_1\odot\cdots \odot v_p, \xi^1\odot\cdots \odot \xi^q\rangle=
\begin{cases} \iota_{v_1}\cdots \iota_{v_p} (\xi^1\odot\cdots \odot \xi^p), & p=q,\\
0, & \mbox{otherwise.}
\end{cases}
$$
The pairing between $S^\bullet(V^\vee)$ and $S^\bullet(V)$ is similarly defined, and we have
\begin{align*}
  \langle v_1\odot\cdots\odot v_n, \xi^1\odot\cdots\odot\xi^n \rangle = (-1)^{(\sum_{i=1}^n\abs{v_i})(\sum_{j=1}^n\abs{\xi^j})}\langle \xi^1\odot\cdots\odot\xi^n, v_1\odot\cdots\odot v_n \rangle.
\end{align*}

Let %us introduce the notation
$$
\O(V) =  \widehat{S}^\bullet(V^\vee)
$$
be the space of formal power series on the graded vector space $V$, %From the point of super geometry, an element in $\O(V)$ is a function on the graded manifold $V$, namely an infinite power series centered at the origin of $V$.
which is a local algebra with the unique maximal ideal
$$
\O^+(V)=\O(V)\odot V^\vee=\widehat{S}^{n\geq 1}(V^\vee)=\prod_{n\geq 1} S^n(V^\vee).
$$
Denote by $r^+:\O(V)\rightarrow \O^+(V)$ the obvious projection.

For all $k \geq 0$, consider the product
\be\label{mu}
 \mu^V_{k+1}: S^k(V) \otimes V \rightarrow S^{k+1}(V), \;\; x \otimes v \mapsto x \odot v,\;\;\forall x \in S^k(V), v \in V.
\ee
The dual map will be denoted by
\be\label{Ik+1}
 I^V_{k+1}: S^{k+1}(V^\vee) \rightarrow S^k(V^\vee) \otimes V^\vee.
\ee
The summation of these $I^V_{k+1}$ defines an operator
\be\label{I}
I^V = \sum_{k \geq 0}I^V_{k+1} \circ r^+: \O(V) \rightarrow \O(V) \otimes V^\vee,
\ee
which is in fact the algebraic de Rham operator of the $\k$-algebra $\O(V)$. It is clear that $I^V$ is an $\O(V)$-derivation valued in the $\O(V)$-bimodule $\O(V) \otimes V^\vee$, i.e., for all $\omega,\omega^\prime \in \O(V)$,
\be\label{derivation of I}
 I^V(\omega \odot \omega^\prime) = \omega \odot I^V(\omega^\prime) + I^V(\omega) \odot \omega^\prime = \omega \odot I^V(\omega^\prime) + (-1)^{\abs{\omega}\abs{\omega^\prime}}\omega^\prime \odot I^V(\omega).%\;\;\forall \omega,\omega^\prime \in \O(V).
\ee
%\begin{Def}
A degree $n$ derivation $D$ of $\O(V)$ is a degree $n$ $\k$-linear map $D:~\O(V) \rightarrow \O(V)$ such that the following Leibniz rule holds:
  \bd
  D(\xi\odot\eta) = D(\xi)\odot\eta + (-1)^{n\abs{\xi}}\xi\odot D(\eta),\;\;\;\forall \xi,\eta \in \O(V).
  \ed
%\end{Def}
The space $\Der(\O(V))$ of derivations of $\O(V)$, together with the graded commutator
\bd
[D_1,D_2] :=  D_1 \circ D_2 - (-1)^{\abs{D_1}\abs{D_2}}D_2 \circ D_1,\;\;\forall D_1,D_2 \in \Der(\O(V)),
\ed
%Note that $[X_1,X_2]$ is a derivation of degree $\abs{X_1}+\abs{X_2}$. Such homogeneous derivations span the set $\Der(\O(V))$, which
is a graded Lie algebra. %In general, a derivation of $\O(V)$ is of the form $\sum \xi_i\otimes \iota_{v_i}$, $\xi_i\in \O(V)$, $v_i\in V$. In other words, $\Der(\O(V))\cong \O(V)\otimes V$. % with commutators as Lie brackets.
%\begin{Ex}
%For any $v \in V$, the left contraction $\iota_v:~\O(V) \rightarrow \O(V)$ is a degree $\abs{v}$ derivation of $\O(V)$. And the space of left contractions is an abelian subalgebra of $\Der(\O(V))$ since $[\iota_u,\iota_v] = 0, \;\; \forall u,v \in V$.
%\end{Ex}

\subsection{Strongly homotopy Lie algebras}
%In this paper, by strongly homotopy Lie algebras, or simply SH Lie algebras, we mean $L_\infty[1]$-algebras defined as follows:
\begin{Def}\label{Def:L[1]-algebra}
  An SH Lie algebra (or $L_\infty[1]$-algebra) is a pair
 $(L,\{\lambda_k\}_{k=0}^\infty)$, simply denoted by $(L,\lambda_\bullet)$, where $L$ is a graded vector space and $\lambda_k:~S^k(L) \rightarrow L, k \geq 0$, called the $k$th-bracket, are degree $1$ linear maps satisfying the following generalized Jacobi identities:
  \be\label{Jacobi1}
    \sum_{k+l = n}\sum_{\sigma \in \sh(l,k)}\epsilon(\sigma)\lambda_{k+1}(\lambda_l(u_{\sigma(1)},\cdots,u_{\sigma(l)}),\cdots,u_{\sigma(n)}) = 0,
  \ee
  for all $k,l,n \geq 0$ and homogeneous elements $u_i \in L, 1 \leq i \leq n$.   Here $\sh(s,k)$ is the set of $(s,k)$-unshuffles.
\end{Def}

%\begin{Rm}
%  The first element $\lambda_0 \in L^1$ is a constant vector. If $n=0$, then Equation \eqref{Jacobi1} becomes $\lambda_1(\lambda_0) = 0$. The generalized Jacobi identities for $n=1, 2$ are
%\begin{align*}
%  n&=1,  &    &\lambda_1(\lambda_1(u)) + \lambda_2(\lambda_0,u) = 0; \\
%  n&=2,  &    &\lambda_3(\lambda_0,u_1,u_2) + \lambda_2(\lambda_1(u_1),u_2) + (-1)^{\abs{u_1}\abs{u_2}}\lambda_2(\lambda_1(u_2),u_1) + \lambda_1(\lambda_2(u_1,u_2)) = 0.
%\end{align*}

%  Thus, $\lambda_0 \in L^1$ is the obstruction for $\lambda_1$ being of square zero, which is also called the curvature of the $L_\infty[1]$-algebra $L$. We call $L$ flat if $\lambda_0$ vanishes.
%\end{Rm}

\begin{Ex}\label{Endomorphism space}
  Let $E$ be a graded vector space. Then $\End(E)$ together with the graded commutator $[-,-]$ is a Lie  algebra. Thus $\End(E)[1]$ with the shifted commutator $\{-,-\}$:
  \be\label{shifted graded commutator}
  \{\bar{\phi},\bar{\psi}\} := (-1)^{\abs{\phi}}[\phi,\psi] = (-1)^{\abs{\phi}}(\phi \circ \psi - (-1)^{\abs{\phi}\abs{\psi}}\psi \circ \phi),\;\;\forall \bar{\phi},\bar{\psi} \in \End(E)[1],
  \ee
    is an $L_\infty[1]$-algebra with only one nontrivial bracket $\lambda_2 = \{-,-\}$. Here $\phi$ and $\bar{\phi}$ are the same element with different degrees: $\abs{\bar{\phi}}=\abs{\phi}-1$.  {The relation between $\bar{\psi}$ and $\psi$ is similar.}
\end{Ex}

\begin{Rm}
  Our Definition~\ref{Def:L[1]-algebra} of $L_\infty[1]$-algebras is not the more commonly known notion of $L_\infty$-algebras (e.g., see standard texts \cites{LM,LS}). In particular, one should notice the different convention of degrees and  signs of $L_\infty$ and $L_\infty[1]$-algebras. What we adopt is similar to that in \cite{Vor}, where  $\Z_2$ grading is used. For a passage connecting our definition to that in \cites{LM,LS}, we refer to \cite{Vor}*{Remark 2.1}.
\end{Rm}
%\subsection{SH Lie algebras as Q-manifolds}
SH Lie algebras could also be characterized as $Q$-manifolds~\cite{AKSZ}:
\begin{Def}
  A homological vector field on $L$ is a degree $1$ derivation $Q$ on $\O(L)$ such that  $Q^2 = \frac{1}{2}[Q,Q] = 0$.
\end{Def}
\begin{Pro}\label{L-infinity as Q}
  Let $L$ be a graded vector space. Then there is a one-to-one correspondence between $L_\infty[1]$-algebra structures on $L$ and homological vector fields on $L$.
\end{Pro}
In fact, on the one hand, if $(L,\{\lambda_k\}_{k \geq 0})$ is an $L_\infty[1]$-algebra, then we can construct a homological vector field $Q_L$ as follows:

For each $k \geq 0$, the dual of $\lambda_k:~S^k(L)\rightarrow L$ is a map
$\lambda_k^\vee:~L^\vee \rightarrow S^k(L^\vee)$, which can be uniquely extended to a degree $1$ derivation $\O(L)\rightarrow \O(L)$. Then define $Q_L$ by
\bd
Q_L = \sum_{k \geq 0}(-1)^k\lambda_k^\vee:~\O(L) \rightarrow \O(L),
\ed
i.e., for all $\xi \in L^\vee, u_i \in L, 1 \leq i \leq k$, we have
\be \label{Q_L}
\langle \xi,\lambda_k(u_1,\cdots,u_k)\rangle=(-1)^{\abs{\xi}+k}\langle Q_L(\xi),u_1\odot \cdots \odot u_k\rangle.
\ee
On the other hand, given a homological vector field $Q_L$ on $L$, we can define a collection of degree $1$ linear maps $\lambda_k:~S^k(L) \rightarrow L, k \geq 0$ by
\be\label{derived brackets}
\lambda_k(u_1,\cdots,u_k) = \iota^{-1}\left([[\cdots[[Q,\iota_{u_1}],\iota_{u_2}]\cdots],\iota_{u_k}]\right) \in L,\;\;\forall u_i \in L.
\ee
Here the map
\bd
\iota^{-1}:~\Der(\O(L)) \rightarrow L
\ed
is defined by
\bd
\langle \iota^{-1}(D),\xi \rangle = pr_0 \circ D(\xi) \in \k,\;\;\;\forall D \in \Der(\O(L)), \xi \in L^\vee,
\ed
where
\bd
pr_0:~\O(L) =  \widehat{S}{^\bullet(L^\vee)} \rightarrow S^0(L^\vee) = \k
\ed
is the obvious projection.
It is clear that the map $\iota^{-1}$ is the left inverse of the contraction operator $\iota_{-}:L\rightarrow \Der(\O(L))$ in the sense that $\iota^{-1}(\iota_u)=u$, for all $u\in L$.

The fact that $Q_L$ in Equation~\eqref{Q_L} is of square zero and the fact that $\{\lambda_k\}_{k\geq 0}$ in Equation~\eqref{derived brackets} defines an $L_\infty[1]$-algebra structure on $L$ are equivalent. Details can be found in \cites{Vor,Vor2}.

For this reason, an $L_\infty[1]$-algebra can be denoted by any of the notations $(L,\lambda_\bullet)$,   $(L,Q_L)$ or $(L,\lambda_\bullet\sim Q_L)$.

Now we recall morphisms of SH Lie algebras:
\begin{Def}\label{morphism of SH Lie algebras}
  Let $(L,\lambda_\bullet\sim Q_L)$ and $(L^\prime, \lambda_\bullet^\prime\sim Q_{L^\prime})$ be two $L_\infty[1]$-algebras.  An $L_\infty[1]$-morphism from $L$ to $L^\prime$ is a morphism $\phi: \O(L^\prime) \rightarrow \O(L)$ of $\k$-algebras such that
  \be\label{Q-morphism}
  \phi \circ Q_{L^\prime} = Q_L \circ \phi: \O(L^\prime) \rightarrow \O(L).
  \ee
  Equivalently, an $L_\infty[1]$-morphism from $L$ to $L^\prime$ is a family of degree zero linear maps
  \bd
  f_k:~S^k(L) \rightarrow L^\prime, \;\; k \geq 0
  \ed
  satisfying the following two conditions:
  \begin{enumerate}[$(1)$]
    \item The element $f_0 \in (L^\prime)^0$ satisfies
    \be \label{f_0-relation}
     \sum_{k \geq 1}\frac{1}{k!}\lambda_k^\prime(f_0,\cdots,f_0) + \lambda_0^\prime = f_1(\lambda_0).
    \ee
    \item For each $n \geq 1$, the relation
    \begin{align}\label{morphism relation}
    &\quad \sum_{k+l=n}\sum_{\sigma \in \sh(l,k)} \epsilon(\sigma)f_{k+1}(\lambda_l(u_{\sigma(1)},\cdots,u_{\sigma(l)})\cdots,u_{\sigma(n)}) \\
    &=\sum_{\substack{i_1,\cdots,i_r \geq 1 \\ i_1 + \cdots + i_r = n}}\sum_{\tau \in \sh(i_1,\cdots,i_r)}\sum_{j \geq 0} \frac{1}{(r+j)!}\epsilon(\tau)\lambda_{r+j}^\prime(f_0^{\odot j},f_{i_1}(u_{\tau(1)},\cdots,u_{\tau(i_1)}),\cdots f_{i_r}(\cdots, u_{\tau(n)}))\notag
    \end{align}
    holds, where $k,l \geq 0$ and $u_i \in L$ are homogeneous.
  \end{enumerate}
\end{Def}
For completeness, we will give a proof on the equivalence of the two definitions of morphisms in Appendix~\ref{appendix}.

\subsection{Connections, curvatures and modules of SH Lie algebras}
This part can be thought of as formal differential geometry of $L_\infty[1]$-algebras. What we shall deal with, namely connections and curvatures, are defined in the same manner of super-connections in \cite{Quillen}.  Some closely related contents can be found in \cites{AC,LM,Luca,Mehta}.

Let us fix an $L_\infty[1]$-algebra $(L,\lambda_\bullet\sim Q_L)$ and a graded vector space $E$.
%The space $\O(L)\otimes E$ can be seen as the space of $E$-valued functions on $L$, which is obviously an $\O(L)$-bimodule: Any element $\omega \in \O(L)$ acts on $\xi \otimes e \in \O(L)\otimes E$ from the left and from the right by
%\begin{align*}
% \omega \odot (\xi \otimes e) &:= (\omega \odot \xi) \otimes e \in \O(L) \otimes E, & (\xi \otimes e) \odot \omega &:= (-1)^{\abs{e}\abs{\omega}}(\xi \odot \omega) \otimes e \in \O(L) \otimes E,
%\end{align*}
%respectively.
A degree $n$ operator $\partial:~\O(L)\otimes E\rightarrow \O(L)\otimes E$ is called an $E$-derivation if it is $\k$-linear and there exists a degree $n$ derivation $\underline{\partial}\in \Der(\O(L))$, called the symbol of $\partial$, such that
$$
\partial(\omega\otimes e)=\underline{\partial}(\omega)\otimes e+(-1)^{n|\omega|}\omega\odot \partial(e),\quad \quad \omega\in \O(L),e\in E.
$$
%The element $\underline{\partial}$ is called the symbol of $\partial$.
Let us denote the space of all $E$-derivations by $\Der(\O(L)\otimes E)$. There is a Lie bracket on $\Der(\O(L)\otimes E)$ defined by its graded commutator
$$
[\partial,\partial^\prime]=\partial\circ \partial^\prime -(-1)^{|\partial||\partial^\prime|}\partial^\prime \circ \partial,\;\forall \partial, \partial^\prime \in \Der(\O(L) \otimes E).
$$
%Meanwhile, it is easy to see that the symbol map
%\bd
% s: \Der(\O(L) \otimes E) \rightarrow \Der(\O(L)),\;\; s(\partial) = \underline{\partial}
%\ed
%is a morphism of Lie algebras. And the kernel of $s$ is the subspace in $\Der(\O(L))\otimes E)$ consisting of $E$-derivations with zero symbols, which can be identified with $\O(L)\otimes \End(E)$, the space of $\O(L)$-linear endomorphisms of $\O(L)\otimes E$. It is also a Lie algebra with the Lie bracket defined by
%$$
% [\omega\otimes \phi, \omega^\prime \otimes \phi^\prime]:= \omega \odot \omega^\prime \otimes [\phi,\phi^\prime],\;\forall \omega\otimes \phi, \omega^\prime \otimes \phi^\prime \in \O(L) \otimes \End(E)
%$$
%where $[\phi,\phi^\prime]$ is the graded commutator in $\End(E)$.

\begin{Lem}\label{Lem:E-derivation}
The Lie algebra $\Der(\O(L) \otimes E)$ is isomorphic to the semidirect product of $\Der(\O(L))$ and $\O(L) \otimes \End(E)$.
\end{Lem}
\bp
Note that $\O(L) \otimes \End(E)$ consists of $E$-derivations with zero symbol, thus a Lie subalgebra of $\Der(\O(L) \otimes E)$. And $\Der(\O(L))$ is also a Lie subalgebra by the inclusion
\bd
 j: \Der(\O(L)) \rightarrow \Der(\O(L) \otimes E), \quad j(X)(\omega\otimes e) = X(\omega) \otimes e,\;\forall X \in \Der(\O(L)), \omega \in \O(L), e \in E,
\ed
which gives rise to a natural splitting of the short exact sequence of Lie algebras
\bd%\label{SES of lie algebras}
 0 \rightarrow \O(L) \otimes \End(E) \hookrightarrow \Der(\O(L) \otimes E) \rightarrow \Der(\O(L)) \rightarrow 0,
\ed
where the Lie algebra morphism $\Der(\O(L) \otimes E) \rightarrow \Der(\O(L))$ is taking symbols. %Moreover, $\Der(\O(L))$ acts on $\O(L) \otimes \End(E)$ by
%\bd
% X(\omega \otimes \phi) = X(\omega) \otimes \phi,\;\;\forall X \in \Der(\O(L)), \omega \in \O(L), \phi \in \End(E).
%\ed
%It can be easily verified that these operations
It gives rise to a semidirect product $\Der(\O(L)) \ltimes (\O(L) \otimes \End(E))$ that is isomorphic to $\Der(\O(L) \otimes E)$.
\ep

\subsubsection{Connections and curvatures}
\begin{Def}
  An $L$-connection on $E$ is a degree $1$ $E$-derivation $\nabla\in \Der(\O(L)\otimes E)$ whose symbol is $Q_L$, i.e., the following Leibniz rule holds:
  \bd
  %\label{Leibniz rule in connection}
  \nabla(\omega \otimes e) = Q_L(\omega) \otimes e + (-1)^{\abs{\omega}}\omega \odot \nabla(e),\;\;\forall \omega \in \O(L), e \in E.
  \ed
  The degree $2$ $E$-derivation
  \bd
  R^\nabla := \nabla^2=\frac{1}{2}[\nabla,\nabla]:~\quad \O(L) \otimes E \rightarrow \O(L) \otimes E
  \ed
  is of zero symbol, i.e., $\O(L)$-linear, and will be called the curvature of $\nabla$. An $L$-connection $\nabla$ is said to be flat if its curvature $R^\nabla$ vanishes.
\end{Def}

%\begin{prop-def}
%The degree $2$ $E$-derivation
%  \bd
%  R^\nabla := \nabla^2=\frac{1}{2}[\nabla,\nabla]:~\quad \O(L) \otimes E \rightarrow \O(L) \otimes E
%  \ed
%is of zero symbol, i.e., $\O(L)$-linear, which will be called the \emph{curvature} of $\nabla$. An $L$-connection $\nabla$ is said to be \emph{flat} if its
%curvature $R^\nabla$ vanishes.
%\end{prop-def}
%\bp
%By Lemma~\ref{Lem:E-derivation}, the symbol map $\Der(\O(L) \otimes E) \rightarrow \Der(\O(L))$ is a Lie algebra morphism. Thus the symbol of $[\nabla,\nabla]$ is $[Q_L,Q_L]$, which vanishes since $Q_L$ is a homological vector field by assumption. Therefore, $R^\nabla$ is $\O(L)$-linear, as desired.
%\ep
The difference of two connections is a degree $1$ endomorphism of the $\O(L)$-module  $\O(L)\otimes E$. Thus the set of all $L$-connections on $E$ is an affine space over $(\O(L)\otimes \End(E))^1$.

According to Lemma \ref{Lem:E-derivation}, any $L$-connection $\nabla$ is determined by an element $D^E\in (\O(L)\otimes \End(E))^1$ so that
$\nabla=Q_L+D^E$. An easy computation shows that the curvature has the form
\begin{align*}
  R^\nabla=\nabla^2 &=(Q_L+D^E) \circ (Q_L+D^E)%\\ &=Q_L^2+D^E \circ Q_L+Q_L\circ D^E+(D^E)^2\\ &
  =Q_L(D^E)+(D^E)^2.
\end{align*}

\begin{Lem}
We have the Bianchi identity:
 \be\label{Bianchi identity}
     Q_L(R^\nabla) + [D^E,R^\nabla] = 0.
  \ee
\end{Lem}
\bp
This follows from straightforward computations:
\bd
[Q_L + D^E, R^\nabla] = [\nabla,\nabla^2]=\nabla\circ \nabla^2-(-1)^{1\times 2}\nabla^2\circ \nabla = 0.
\ed
\ep

%Given any $\partial\in \Der(\O(L)\otimes E)$, it induces the dual $E^\vee$-derivation $\partial^\vee$, with the same symbol $\underline{\partial}$, such that
%$$
%\langle \partial^\vee(\beta),\alpha\rangle=\underline{\partial}\langle \beta,\alpha\rangle-
%(-1)^{\abs{\partial}\abs{\beta}}\langle \beta, \partial \alpha\rangle, \quad \quad \forall \alpha\in \O(L)\otimes E, \beta\in \O(L)\otimes E^\vee.
%$$
%In particular,
Given an $L$-connection $\nabla$ on $E$, there corresponds a dual $L$-connection $\nabla^\vee:~\O(L) \otimes E^\vee \rightarrow \O(L) \otimes E^\vee$ on $E^\vee$. Explicitly, we have
\be\label{dual connection}
 \langle \nabla^\vee(f), g \rangle = Q_L\langle f,g \rangle - (-1)^{\abs{f}} \langle f, \nabla(g) \rangle \in \O(L), \quad \quad \forall f \in \O(L) \otimes E, g \in \O(L) \otimes E^\vee.
\ee

For two graded vector spaces with $L$-connections $(E,\nabla^E)$ and $(F, \nabla^F)$, the induced $L$-connection on $E \otimes F$ is given by
\begin{equation}\label{Eqt:nablaproduct}
\nabla^{E\otimes F}(\omega\otimes e\otimes f)=Q_L(\omega)\otimes e\otimes f+(-1)^{\abs{\omega}}\omega \odot (\nabla^E(e) \otimes f) +(-1)^{\abs{\omega}+\abs{e}}(\omega\otimes e)\otimes_{\O(L)} \nabla^{F}(f).
\end{equation}
Here we used the canonical isomorphism
\begin{align*}
  (\O(L) \otimes E) \otimes_{\O(L)} (\O(L) \otimes F) &\xrightarrow{\cong} \O(L) \otimes (E \otimes F)% \\
     %      (\omega \otimes e) \otimes_{\O(L)} (\omega^\prime \otimes e^\prime) &\mapsto (-1)^{\abs{e}\abs{\omega^\prime}}\omega \odot \omega^\prime \otimes (e \otimes e^\prime)
\end{align*}
to view the last term $(\omega\otimes e)\otimes_{\O(L)} \nabla^{F}(f)$ as an element in $\O(L) \otimes (E \otimes F)$.

%In particular, an $L$-connection $\nabla$ on $E$ determines an $L$-connection $\nabla^\vee$ on $E^\vee$, as well as an $L$-connection $\nabla^{\End(E)}$ on $\End(E)\cong E^\vee \otimes E \cong E\otimes E^\vee$.

%\begin{Lem}\label{Lem:nablaEndE} We have
%$$\nabla^{\End(E)}=[\nabla,-],$$ as a map $\O(L)\otimes \End(E)\rightarrow \O(L)\otimes \End(E)$.
%\end{Lem}
%\bp
%Taking any element $\omega\otimes e^\vee \otimes f\in \O(L)\otimes \End(E)\cong \O(L)\otimes E^\vee \otimes E$ and $e\in E$, we have
%\begin{align*}
%  &\quad \langle \nabla^{\End(E)}(\omega\otimes e^\vee \otimes f), e \rangle\\
%  &= \langle Q_L(\omega)\otimes e^\vee \otimes f+(-1)^{\abs{\omega}}\omega\otimes  D^{E^\vee}(e^\vee) \otimes f + (-1)^{\abs{\omega}+\abs{e^\vee}}\omega\otimes e^\vee \otimes D^E(f), e \rangle \\
%  &=Q_L(\omega)\otimes \langle e^\vee\otimes f,e\rangle +(-1)^{\abs{\omega}+\abs{e^\vee}}\langle \omega\otimes e^\vee \otimes D^E(f), e\rangle-(-1)^{\abs{\omega}+\abs{f}+\abs{e^\vee}}\langle \omega\otimes e^\vee\otimes f,D^E(e) \rangle \\
%  &=\nabla((\omega\otimes e^\vee\otimes f)(e))-(-1)^{\abs{\omega}+\abs{f}+\abs{e^\vee}}(\omega\otimes e^\vee\otimes f)(D^E(e))\\
%  &=[\nabla, \omega\otimes e^\vee \otimes f](e).
%\end{align*}
%This completes the proof.
%\ep
In particular, we have the induced $L$-connection $\nabla^{\Hom(E,F)}$ on $\Hom(E,F)\cong E^\vee\otimes F$.
%The following lemma follows from straightforward computations:%can be proved in an analogous manner to that of Lemma \ref{Lem:nablaEndE}.
\begin{Lem}\label{Lem:nablaHomEF}
For each $\Psi \in \O(L) \otimes \Hom(E,F) \cong \Hom_{\O(L)}(\O(L) \otimes E, \O(L) \otimes F)$,
\begin{equation}\label{Eqt:nablaHom}
\nabla^{\Hom(E,F)}(\Psi)= \nabla^F \circ \Psi  - (-1)^{\abs{\Psi}} \Psi\circ \nabla^E.
\end{equation}
\end{Lem}

\subsubsection{Modules over SH Lie algebras}
\begin{Def}
 An $L$-module $(E,\partial^E_L)$ is a graded vector space $E$ together with a flat $L$-connection $\partial_L^E:~\O(L)\otimes E\rightarrow \O(L)\otimes E$, which will be called the Chevalley-Eilenberg differential of $E$. The associated cohomology $\Ha^\bullet(\O(L)\otimes E,\partial_L^E)$ will be denoted by $\Ha^\bullet_{\mathrm{CE}}(L,E)$.
\end{Def}
%The operator $\partial_L^E$ generalizes the well-known Chevalley-Eilenberg differential. For this reason, we will also call $\partial_L^E$ the \emph{Chevalley-Eilenberg differential} (simply, the differential), and $(\O(L)\otimes E,\partial_L^E)$ the \emph{Chevalley-Eilenberg complex of $L$ with coefficient $E$}. Accordingly, the cohomology space $\Ha^\bullet(\O(L)\otimes E,\partial_L^E)$ is denoted by $\Ha^\bullet_{\mathrm{CE}}(L,E)$.
Recall that as an $L$-connection, $\partial_L^E$ is determined by an element $D^E \in (\O(L)\otimes \End(E))^1$ such that $\partial_L^E=Q_L+D^E$. The flat condition of $\partial_L^E$ becomes a Maurer-Cartan equation:
\be\label{MC-EQ}
 Q_L(D^E)+(D^E)^2=0.
\ee
We have alternative descriptions of $L$-modules.
\begin{Pro}\label{Thm:module}
  Let $L = (L,\lambda_\bullet\sim Q_L)$ be an $L_\infty[1]$-algebra, and $E$ be a graded vector space. The following data are mutually determined:
  \begin{enumerate}[$(1)$]
    \item An $L$-module structure on $E$;
    \item An $L$-module structure on $E^\vee$;
    \item A degree $1$ derivation $Q\in \Der(\O(L\oplus E))$ such that $Q|_{\O(L)}=Q_L$, $Q(E^\vee)\subset \O(L) \otimes E^\vee$ and $Q^2=0$;
    \item An abelian extension $L \oplus E$ of $L$ along $E$ ({also called a} \emph{semi-direct product} of $L$ with $E$):~ i.e., $L \oplus E$ carries an $L_\infty[1]$-algebra structure $\{\tilde{\lambda}_k\}_{k \geq 0}$ such that (i) $L$ is an $L_\infty[1]$ subalgebra (in particular,   $\tilde{\lambda}_0={\lambda}_0 \in   L^1$); (ii) $E$ is an ideal:~$\tilde{\lambda}_k(E,\cdots)\subset E, k\geq 1$; (iii) $E$ is abelian:~$\tilde{\lambda}_k(E,E,\cdots)=0, k\geq 2$;
    \item A family of degree $1$ linear maps
    \be\label{m_k}
     m^E_k:~S^{k-1}(L) \otimes E \rightarrow E,
    \ee
    $k=1,2,\cdots$, such that
     \begin{align}\label{module relation}
    &\sum_{k+l=n}\sum_{\sigma \in \sh(l,k)}\epsilon(\sigma)m^E_{k+2}(\lambda_l(u_{\sigma(1)},\cdots,u_{\sigma(l)}),\cdots,u_{\sigma(n)},e) \notag\\
    &=-\sum_{k+l=n}\sum_{\tau \in \sh(k,l)}\epsilon(\tau)(-1)^{\dagger^\tau_k} m^E_{k+1}(u_{\tau(1)},\cdots,u_{\tau(k)},m^E_{l+1}(u_{\tau(k+1)}, \cdots, u_{\tau(n)},e))
    \end{align}
     hold for all $k,l,n \geq 0$ and homogeneous elements $u_i \in L, e \in E$, where $\dagger_k^\tau = \sum_{i=1}^k\abs{u_{\tau(i)}}$;
    \item An $L_\infty[1]$-morphism $f^E = \{f^E_k\}$ from $L$ to $\End(E)[1]$, where
    \be\label{modulebyf}
     f^E_{k}: S^k(L)\rightarrow \End(E)[1], \qquad\quad k \geq 0.
    \ee
    \end{enumerate}
\end{Pro}
\bp
For the equivalence $(1) \Leftrightarrow (2)$, we need to show that an $L$-connection $\partial_L^E$ on $E$ is flat if and only if its dual $L$-connection $\partial_L^{E^\vee}$ on $E^\vee$ is flat.
In fact, by Equation~\eqref{dual connection}, we have
\begin{align*}
  \langle (\partial_L^{E^\vee})^2(e^\vee), e \rangle &= Q_L \langle  \partial_L^{E^\vee}(e^\vee), e \rangle - (-1)^{\abs{e^\vee}+1}\langle  \partial_L^{E^\vee}(e^\vee), \partial_L^E(e) \rangle \\
  &= Q_L \langle \partial_L^{E^\vee}(e^\vee), e \rangle - (-1)^{\abs{e^\vee}+1}Q_L\langle e^\vee, \partial_L^E(e) \rangle -\langle e^\vee, (\partial_L^E)^2 (e) \rangle \\
  &=-\langle e^\vee, (\partial_L^E)^2 (e) \rangle.
\end{align*}
Thus $(\partial_L^E)^2 = 0$ is equivalent to $(\partial_L^{E^\vee})^2 = 0$.

The equivalence $(2) \Leftrightarrow (3)$ is obvious.

To see the equivalence  $(3) \Leftrightarrow (4)$, we note that:~$Q^2=0 \Leftrightarrow L \oplus E$ is an $L_\infty[1]$-algebra (by Proposition \ref{L-infinity as Q}); $Q|_{\O(L)}=Q_L \Leftrightarrow$ (i) in (4); $Q(E^\vee) \subset \O(L) \otimes E^\vee \Leftrightarrow$ (ii) and (iii) in (4).

To see $(4) \Leftrightarrow (5)$, we note that by setting
\bd
m^E_k(u_1,\cdots,u_{k-1},e)=\tilde{\lambda}_k(u_1,\cdots,u_{k-1},e),
\ed
Equation \eqref{module relation} is equivalent to the generalized Jacobi identity \eqref{Jacobi1} of $(L\oplus E, \tilde{\lambda}_\bullet)$.  %when applied to elements $u_1,\cdots,u_{n}\in L,e\in E$. Other situations are merely the generalized Jacobi identity of $L$ itself.

Finally, we show that $(5) \Leftrightarrow (6)$. In fact, for all $k \geq 0$, $f^E_k$ and $m^E_{k+1}$ are mutually determined by
\bd
 f^E_k(u_1,\cdots,u_k)(e) = m^E_{k+1}(u_1,\cdots,u_k,e).
\ed
Now we reformulate Equation~\eqref{module relation} in terms of $f^E_k$: For $n = 0$, it becomes
\bd
m^E_2(\lambda_0,e) = -m^E_1(m^E_1(e)),\;\;\forall e \in E,
\ed
which is equivalent to
\bd
f^E_1(\lambda_0)(e) = -\frac{1}{2}[f^E_0,f^E_0](e) = \frac{1}{2}\{f^E_0,f^E_0\}(e), \;\;\forall e \in E.
\ed
Here $[-,-]$ is the graded commutator in $\End(E)$ and $\{-,-\}$ is the shifted graded commutator in $\End(E)[1]$ (see Example~\ref{Endomorphism space}), which is the only nontrivial bracket in $\End(E)[1]$.

For $n \geq 1$, Equation \eqref{module relation} can be reorganized as
\begin{align}
  &\sum_{k+l=n}\sum_{\sigma \in \sh(l,k)}\epsilon(\sigma)f^E_{k+1}(\lambda_l(u_{\sigma(1)},\cdots,u_{\sigma(l)}),\cdots,u_{\sigma(n)}) \notag\\
  \label{equationonf} &= \frac{1}{2}\{f^E_0,f^E_n(u_1,\cdots,u_n)\} + \sum_{k+l=n}\sum_{\tau \in \sh(l,k)}\epsilon(\tau)\frac{1}{2}\{f^E_k(u_{\tau(1)},\cdots,u_{\tau(k)}), f^E_l(u_{\tau(k+1)},\cdots,u_{\tau(n)})\}.
\end{align}
These are exactly Equation~\eqref{morphism relation}.
\ep

\begin{Rm}
According to this theorem, we can write the explicit relation between $D^E$ (or $\partial_L^E$) and $m^E_\bullet$:
\be\label{temp6}
m^E_{k+1}(u_1,\cdots,u_k,e)=(-1)^{\abs{e}\ast_k+k+1}D^E(e)\llcorner u_1\llcorner\cdots \llcorner u_k,\; \forall u_1,\cdots,u_k\in L,e\in E,
\ee
where $\ast_k=\sum_{i=1}^k \abs{u_i}$.
\end{Rm}

From now on, we will denote an $L$-module by any of the notations $(E,\partial_L^E)$, $(E,D^E)$, $(E, m^E_\bullet)$, $(E,f_\bullet^E)$ or $(E,\partial_L^E\sim m^E_\bullet)$, etc. The data $\partial_L^E:\O(L)\otimes E\rightarrow \O(L)\otimes E$, $D^E\in (\O(L)\otimes \End(E))^1$, $\{m^E_k\}_{k\geq 1}$ in Equation~\eqref{m_k}, and $\{f^E_k\}$ in Equation~\eqref{modulebyf} are mutually determined by the above theorem.

\begin{Ex}[Adjoint module] \label{Adjoint rep.}
  Let $(L,\lambda_\bullet)$ be an $L_\infty[1]$-algebra. Then the maps
  \be\label{adjoint m}
  m^L_{k+1} = \lambda_{k+1} \circ \mu^L_{k+1}: S^k(L) \otimes L \rightarrow L, \quad k \geq 0
  \ee
  make $L$ an $L$-module, where $\mu^L_{k+1}$ are defined by Equation~\eqref{mu}. In analogy to Lie algebras, we call it the adjoint $L$-module or the adjoint representation of $L$. %Meanwhile, we will also denote by $\partial_L^L = Q_L + D^L$ the corresponding flat $L$-connection on $L$.
\end{Ex}

Now we introduce morphisms of $L$-modules.
%An object in $\Mod_L$ is a pair $(E,\partial_L^E)$, where $E$ is a graded vector space and $\partial_L^E:\O(L)\otimes E\rightarrow \O(L)\otimes E$ is a flat $L$-connection on $E$, namely the Chevalley-Eilenberg differential.

%Given objects $(E,\partial_L^E),(F,\partial_L^F) \in \Mod_L$, the dual $(E^\vee, \partial_L^{E^\vee})$, the product $(E \otimes F, \partial_L^{E\otimes F})$, and the Hom $ (\Hom(E,F) \cong E^\vee \otimes F, \partial_L^{\Hom(E,F)})$ are all objects in $\Mod_L$. The associated differentials are determined by Equations \eqref{dual connection}, \eqref{Eqt:nablaproduct}, \eqref{Eqt:nablaHom}, respectively.

\begin{Def}[Morphisms of $L$-modules]\label{morphisms of modules}
  Let $(E,\partial_L^E \sim m^E_\bullet)$ and $(F,\partial_L^F \sim m^F_\bullet)$ be two $L$-modules. A morphism of $L$-modules from $E$ to $F$ is an $\O(L)$-linear map $\phi:~ \O(L) \otimes E  \rightarrow  \O(L) \otimes F$ of degree $0$ which is also a chain map, i.e.,
  \be\label{equ. de morphism de modules}
  \phi \circ \partial_L^E = \partial_L^F \circ \phi.
  \ee
  Equivalently, an $L$-module morphism from $E$ to $F$ is an element $\phi = \sum_{k \geq 0}\phi_k\in \O(L)\otimes \Hom(E,F)$, where $\phi_k \in S^k(L^\vee) \otimes  {\Hom}(E,F)$, such that,
  for all $n \geq 0$,
  \begin{align}
  &\sum_{k+l=n}\sum_{\sigma \in \sh(k,l)}\epsilon(\sigma)(-1)^{\ast_n(\dagger^\sigma_k+1)+k}\phi_l(u_{\sigma(k+1)}, \cdots,u_{\sigma(n)})(m^E_{k+1}(u_{\sigma(1)},\cdots,u_{\sigma(k)},e)) \notag\\
  &= \sum_{k+l=n}\sum_{\tau \in \sh(k,l)}\epsilon(\tau)(-1)^{k}\phi_{n-k+1}(\lambda_k(u_{\tau(1)},\cdots,u_{\tau(k)}),\cdots,u_{\tau(n)})(e) \notag \\
  &\quad + \sum_{k+l=n}\sum_{\tau \in \sh(k,l)}\epsilon(\tau)(-1)^{(\ast_n+1)\dagger^\tau_k+n-k}m^F_{n-k+1}(u_{\tau(k+1)}, \cdots, u_{\tau(n)},\phi_k(u_{\tau(1)},\cdots,u_{\tau(k)})(e)) \label{module morphism eq}
 \end{align}
 holds for all $e \in E, u_i \in L, k,l \geq 0$, where $\ast_n = \sum_{i=1}^n\abs{u_i}, \dagger_k^\sigma = \sum_{i=1}^k\abs{u_{\sigma(i)}}$ and $\dagger_k^\tau = \sum_{i=1}^k\abs{u_{\tau(i)}}$.

  The set of such morphisms will be denoted by $\Hom_L(E,F)$.
\end{Def}

Let $(E,f_\bullet^E)$ be an $L$-module. Consider the family of maps
\be\label{moduleasmorphism}
 \phi^{L,E}_k = f^E_{k+1} \circ \mu_{k+1}^L: S^k(L) \otimes L \rightarrow \End(E)[1],\;\;\; k \geq 0,
\ee
where $\mu_{k+1}^L: S^k(L) \otimes L \rightarrow S^{k+1}(L)$ is defined by Equation~\eqref{mu}. It turns out to be a canonical morphism of $L$-modules:
\begin{Lem}
The family of maps $\phi^{L,E} = \{\phi^{L,E}_\bullet\}$ defines an $L$-module morphism from the adjoint $L$-module $L$ to $\End(E)[1]$.
\end{Lem}
\bp
 To see that $\phi^{L,E}$ is an $L$-module morphism, it suffices to show that $\phi^{L,E}$ satisfies Equation~\eqref{module morphism eq}. In fact, note that
 \bd
  m^{\End(E)[1]}(\varphi) = \{m^E,\varphi\},\;\;\forall \varphi \in \End(E)[1],
 \ed
 where $\{-,-\}$ is the shifted graded commutator~\eqref{shifted graded commutator} in $\End(E)[1]$. Then Equation~\eqref{module morphism eq} follows by reformulating Equation~\eqref{equationonf} via Equations~\eqref{adjoint m} and~\eqref{moduleasmorphism}.
\ep

\section{Atiyah classes of SH Lie pairs}\label{Atiyah classes}
\subsection{SH Lie pairs}
\begin{Def}
  By an SH Lie pair $(L,A)$, we mean an $L_\infty[1]$-algebra  $(L,\lambda_\bullet\sim Q_L)$ with a subalgebra $A \subset L$. The structure maps in $A$ are again denoted by $\{\lambda_k\}_{k \geq 0} $. In particular, $\lambda_0 \in A^1 \subset L^1$.
\end{Def}
The homological vector field $Q_A: \O(A)\rightarrow \O(A)$ on the subalgebra $A$ is determined by $Q_L$. In fact, the inclusion map $j:A\rightarrow L$ gives rise to $j^\vee:L^\vee\rightarrow A^\vee$ and a surjective morphism of commutative algebras $j^\vee:\O(L)\rightarrow \O(A)$. The condition that $A$ is a subalgebra in $L$ is equivalent to
\bd
j^\vee \circ Q_L = Q_A \circ j^\vee:~\O(L) \rightarrow \O(A).
\ed
We fix such an  SH Lie pair $(L,A)$. For simplicity, here and in the sequel, we write $B=L/A$.
\begin{Lem}\label{Lem:AmodonB}
  The quotient space $B$ is a canonical $A$-module.
\end{Lem}
\bp
There is an exact sequence of graded vector spaces
\be\label{SES de GVS}
\xymatrix{
0 \ar[r] & A \ar[r]^-{j} & L \ar[r]^-{p} & B \ar[r] & 0.
}
\ee
The canonical $A$-module structure on $B$
\bd
m_k:~S^{k-1}(A) \otimes B \rightarrow B, k \geq 1
\ed
is defined by
\bd
m_k(a_1,\cdots,a_{k-1},b) = p \circ \lambda_k(a_1,\cdots,a_{k-1},l), \;\;\forall a_i \in A, l\in L \mbox{~such that~} p(l)=b.
\ed
These $m_k, k\geq 1$, are well-defined because $A$ is a subalgebra.
That $\{m_k\}_{k \geq 1}$ satisfies Equation~\eqref{module relation} follows from the generalized Jacobi identity \eqref{Jacobi1}.
\ep

It follows that the dual vector space
\bd
B^\vee=(L/A)^\vee\cong A^\perp=\ker (j^\vee:~L^\vee \rightarrow A^\vee)
\ed
is also an $A$-module, which will be denoted by $(A^\perp,\partial^\perp_A)$.

Define an operator $J: \O(L) \rightarrow \O(A) \otimes L^\vee$ by the commutative diagram:
\bd
\xymatrix{
  \O(L) \ar[r]^-{I^L} \ar[dr]_-{J} & \O(L) \otimes L^\vee \ar[d]^-{j^\vee \otimes 1}\\
                                 & \O(A) \otimes L^\vee.
}
\ed
Here $I^L: \O(L) \rightarrow \O(L) \otimes L^\vee$ is defined by Equation~\eqref{I}. It follows immediately that
\be\label{IJ}
(1 \otimes j^\vee) \circ J = j^\vee \circ I^L = I^A \circ j^\vee:~\O(L)\rightarrow \O(A) \otimes A^\vee.
\ee
It is also easy to see that $J$ is a derivation valued in the $\O(L)$-bimodule $\O(A) \otimes L^\vee$, i.e., for all $ \omega,\omega^\prime \in \O(L)$,
\be\label{Jderivation}
 J(\omega \odot \omega^\prime) = j^\vee(\omega) \odot J(\omega^\prime) +  J(\omega) \odot j^\vee(\omega^\prime) = j^\vee(\omega) \odot J(\omega^\prime) + (-1)^{\abs{\omega}\abs{\omega^\prime}}j^\vee(\omega^\prime) \odot J(\omega). %\;\;\forall \omega,\omega^\prime \in \O(L).
\ee
For any graded vector space $E$, $\O(L) \otimes \End(E)$ has an obvious associative product $\circ$. By Equation~\eqref{Jderivation}, the map
\bd
 J \otimes 1: \O(L) \otimes \End(E) \rightarrow \O(A) \otimes L^\vee \otimes \End(E)
\ed
satisfies, for all $\phi,\psi \in \O(L)\otimes\End(E)$,
\be\label{J1}
 (J\otimes 1)(\phi\circ\psi) = (j^\vee\otimes 1)(\phi)\circ (J\otimes 1)(\psi) + (J\otimes 1)(\phi)\circ(j^\vee\otimes 1)(\psi).
\ee

\begin{Pro}\label{dualAmodonBvee}
Let $\partial^\perp_A = Q_A+D^\perp$ be the dual $A$-module structure on $B^\vee \cong A^\perp$. Then for any $\omega \in \ker(j^\vee) \subset \O(L)$, we have $J(\omega) \in \O(A) \otimes A^\perp$ and
\be\label{key relation}
 \partial_A^\perp \left(J(\omega)\right) = J\left(Q_L(\omega)\right).
\ee
\end{Pro}
\bp
Let $(B,D^B \sim m_\bullet)$ be the $A$-module structure on $B$ in Lemma~\ref{Lem:AmodonB}. We first show that
\be\label{Dperp}
D^\perp(\xi) = J\left(Q_L(\xi)\right),\;\;\;\;\forall \xi \in A^\perp.
\ee
In fact, using Equation~\eqref{temp6}, we have, for all $p(l) \in B, a_1\odot\cdots\odot a_k \in S^k(A)$,
\begin{align*}
 &\quad \langle D^\perp(\xi), p(l) \rangle (a_1\odot\cdots\odot a_k) = (-1)^{\abs{\xi}+1}\langle \xi, D^B(p(l))(a_1\odot\cdots\odot a_k)\rangle \\
 &= (-1)^{\abs{\xi}+\abs{l}\ast_k+k}\langle \xi, m_{k+1}(a_1,\cdots,a_k,p(l)) = (-1)^{\abs{\xi}+\abs{l}\ast_k+k}\langle \xi, p \circ \lambda_{k+1}(j(a_1),\cdots,j(a_k),l) \rangle \\
 &= (-1)^{k+1}\langle (j^\vee \circ I^L \circ \lambda_{k+1}^\vee)(\xi),p(l) \rangle (a_1\odot\cdots\odot a_k) \\
 &= \langle (J \circ Q_L)(\xi), p(l) \rangle(a_1\odot\cdots\odot a_k),
\end{align*}
where $\ast_k = \sum_{i=1}^k\abs{a_i}$, since $Q_L = \sum_{k\geq 0} (-1)^k\lambda_k^\vee$. This proves Equation~\eqref{Dperp}.

To prove Equation~\eqref{key relation}, it suffices to consider elements of the form $\omega \odot \xi \in \O(L) \odot A^\perp \cong \ker(j^\vee)$, where $\omega \in \O(L), \xi \in A^\perp$. Then using Equation~\eqref{Dperp}, we have
%In fact, since $Q_L$ preserves $\ker(j^\vee)$ and $J$ is a derivation valued in the $\O(L)$-bimodule $\O(A) \otimes L^\vee$,
\begin{align*}
   &\quad J \left(Q_L(\omega \odot \xi)\right) = J(Q_L(\omega)\odot\xi+(-1)^{\abs{\omega}}\omega\odot Q_L(\xi)) \\
   &= j^\vee \left(Q_L(\omega)\right)\odot J(\xi)+ J \left(Q_L(\omega)\right)\odot j^\vee(\xi) + (-1)^{\abs{\omega}}\left(J(\omega)\odot j^\vee\left( Q_L(\xi)\right) + j^\vee(\omega)\odot J\left(Q_L(\xi)\right)\right) \\
   &= Q_A \left(j^\vee(\omega)\right)\odot\xi + (-1)^{\abs{\omega}}j^\vee(\omega)\odot D^\perp(\xi) = \partial_A^\perp \left( J(\omega \odot \xi)\right).
\end{align*}
This completes the proof.
\ep

\subsection{Construction of Atiyah classes}
In this part, besides the  SH Lie pair $(L,A)$, we fix an $A$-module $(E,m_\bullet\sim D^{A,E})$.
As usual, assume that the differential is of the form
$$
\partial_A^E = Q_A + D^{A,E},
$$
where $D^{A,E} \in (\O(A) \otimes \End(E))^1$ is the $A$-module structure on $E$, and it will be treated as an $\O(A)$-linear map $\O(A)\otimes E\rightarrow \O(A)\otimes E$.

Meanwhile, $A^\perp \otimes \End(E)$ carries an $A$-module structure,  its differential
\bd
\partial_A^{A^\perp\otimes \End(E)}:~~ \O(A) \otimes (A^\perp \otimes \End(E)) \rightarrow \O(A) \otimes (A^\perp \otimes \End(E))
\ed
is expressed by
\bd
%\label{CE}
\partial_A^{A^\perp\otimes \End(E)}=Q_A+D^\perp+[D^{A,E},-].
\ed
The corresponding cohomology space is denoted by $\Ha_{\mathrm{CE}}^\bullet(A, A^\perp \otimes \End(E))$.

Since $j^\vee: \O(L) \rightarrow \O(A)$ is surjective, one is able to find some $D^{L,E} \in (\O(L) \otimes \End(E))^1$ such that $(j^\vee \otimes 1)(D^{L,E}) = D^{A,E}$. Thus we get an $L$-connection $\nabla = Q_L + D^{L,E}$ on $E$ subject to the commutative diagram
\bd
\xymatrix{
\O(L)\otimes E \ar[r]^-{\nabla} \ar[d]^-{j^\vee\otimes 1} & \O(L)\otimes E \ar[d]^-{j^\vee \otimes 1} \\
\O(A)\otimes E \ar[r]^-{\partial^E_A}                     & \O(A)\otimes E.
}
\ed
We call $\nabla$ an $L$-connection extending $(A,E)$. However, it is not necessarily flat. The curvature of $\nabla$ is easily available:
\bd
R^{\nabla} = Q_L(D^{L,E}) + (D^{L,E})^2 \in \O(L) \otimes \End(E).
\ed

We observe the following commutative diagram:
\bd
\xymatrix{
    E \ar[r]^-{R^\nabla} \ar@{=}[d] & \O(L)\otimes E \ar[r]^-{J\otimes 1}  \ar[d]^-{j^\vee\otimes 1} & \O(A)\otimes L^\vee \otimes E \ar[d]^-{1\otimes j^\vee\otimes 1}\\
    E \ar[r]^-{(\partial_A^E)^2=0} & \O(A)\otimes E \ar[r]^-{I^A \otimes 1} & \O(A) \otimes A^\vee \otimes E,
}
\ed
which implies that
\bd
(1\otimes j^\vee \otimes 1)(J \otimes 1)(R^{\nabla}) = 0.
\ed
Hence, we get an element
\be\label{Atiyah cocycle}
\alpha_{\nabla}^E  := (J \otimes 1)(R^{\nabla}) \in \O(A) \otimes A^\perp \otimes \End(E)
\ee
of degree $2$.

\begin{thm-def}\label{main prop-def}
  \begin{enumerate}[$(1)$]
    \item The element $\alpha_{\nabla}^E$ defined by Equation~\eqref{Atiyah cocycle} is a cocycle of the Chevalley-Eilenberg complex of $A$ with coefficient in $A^\perp \otimes \End(E)$, which will be called the Atiyah cocycle of the  SH Lie pair $(L,A)$ with respect to the $A$-module $E$ and the $L$-connection $\nabla$ extending   $(A,E)$.
    \item The cohomology class $[\alpha^E]=[\alpha_{\nabla}^E] \in \Ha^2_{\mathrm{CE}}(A, A^\perp \otimes \End(E))$ does not depend on the choice of the $L$-connection $\nabla$   extending $(A,E)$. We call it the Atiyah class of the  SH Lie pair $(L,A)$ with respect to the $A$-module $E$.
    \item For the   canonical $A$-module $L/A$, there associates a canonical Atiyah class
    $$[\alpha^{L/A}]\in \Ha^2_{\mathrm{CE}}(A, A^\perp \otimes \End(L/A))=\Ha^2_{\mathrm{CE}}(A,\Hom(L/A\otimes L/A,L/A)).$$
  \end{enumerate}
\end{thm-def}
Before giving a proof of the above theorem, we prove the following
\begin{Lem}\label{pro:temp1}
If $X \in \O(L) \otimes \End(E)$ satisfying $(j^\vee \otimes 1)(X)=0$, then
\begin{align*}
[D^{A,E},(J\otimes 1)(X)] &= (J\otimes 1)[D^{L,E},X].
\end{align*}
%where $[-,-]$ is the graded commutator Lie bracket on $\O(L) \otimes \End(E)$.
\end{Lem}
\bp
Without lose of generality, we may assume that $X$ is homogeneous. Note that $(j^\vee \otimes 1)(D^{L,E}) = D^{A,E}$. We have
\begin{align*}
 &\quad (J\otimes 1)[D^{L,E},X] = (J \otimes 1)(D^{L,E} \circ X - (-1)^{\abs{X}}X \circ D^{L,E}) \\
 &= (j^\vee \otimes 1)(D^{L,E}) \circ (J \otimes 1)(X) + (J \otimes 1)(D^{L,E}) \circ (j^\vee \otimes 1)(X) \\
 &\quad - (-1)^{\abs{X}}((j^\vee \otimes 1)(X) \circ (J \otimes 1)(D^{L,E}) + (J \otimes 1)(X) \circ (j^\vee \otimes 1)(D^{L,E})) \qquad\quad(\text{by Equation~\eqref{J1}}) \\
 &= D^{A,E} \circ (J \otimes 1)(X) - (-1)^{\abs{X}}(J \otimes 1)(X) \circ D^{A,E} = [D^{A,E},(J \otimes 1)(X)].
\end{align*}
\ep

\bp[Proof of Theorem-Definition \ref{main prop-def}]
$(1)$. Note that $(j^\vee\otimes 1)(R^\nabla)=(\partial^E_A)^2=0$. It follows that $(J \otimes 1)(R^\nabla) \in \O(A) \otimes A^\perp \otimes \End(E)$. Thus
\begin{align*}
  \partial_A^{A^\perp\otimes \End(E)}(\alpha_{\nabla}^E)&= (Q_A+D^\perp+[D^{A,E},-])((J\otimes 1)(R^\nabla))\\
                          &= (Q_A + D^\perp)((J\otimes 1)(R^\nabla)) + [D^{A,E},(J\otimes 1)(R^\nabla)]\\
                          &= (J\otimes 1)(Q_L(R^\nabla)+[D^{L,E},R^\nabla])=0,
\end{align*}
where we have used Equation~\eqref{key relation}, Lemma~\ref{pro:temp1}, and the Bianchi identity~\eqref{Bianchi identity} in the last two steps.

$(2)$. Let $\tilde{\nabla} = Q_L + \tilde{D}^{L,E}$ be another $L$-connection extending $(A,E)$. Then
$$
\phi= \nabla - \tilde{\nabla} = D^{L,E} - \tilde{D}^{L,E} \in (\O(L) \otimes \End(E))^1
$$
satisfies
\begin{align*}
  (j^\vee\otimes 1)(\phi) &= 0, & (J\otimes 1)(\phi) &\in \O(A) \otimes A^\perp \otimes \End(E).
\end{align*}
It follows from Equation~\eqref{J1} that $(J \otimes 1)(\phi^2) = 0$. Therefore, we have
\begin{align*}
  \alpha_{\nabla}^E-\alpha_{\tilde{\nabla}}^E &= (J\otimes 1)(R^\nabla - R^{\tilde{\nabla}}) =(J\otimes 1)\left(Q_L(D^{L,E}) + (D^{L,E})^2-Q_L(\tilde{D}^{L,E}) - (\tilde{D}^{L,E})^2\right) \\
  &= (J\otimes 1)(Q_L (\phi)+[\tilde{D}^{L,E},\phi]+\phi^2 ) = (J\otimes 1)(Q_L (\phi)+[\tilde{D}^{L,E},\phi])\\
  &= (Q_A + D^\perp)((J\otimes 1)(\phi))+ [D^{A,E},(J\otimes 1)(\phi)]\quad\quad(\text{by Equation~\eqref{key relation} and Lemma~\ref{pro:temp1}})\\
  &= (Q_A + D^\perp + [D^{A,E},-])((J\otimes 1)(\phi)) = \partial_A^{A^\perp\otimes \End(E)} ((J\otimes 1)(\phi)),
\end{align*}
which implies that $[\alpha_{\nabla}^E] = [\alpha_{\tilde{\nabla}}^E]$.

Finally, statement $(3)$ follows from the standard identification $A^\perp \cong (L/A)^\vee$.
\ep

We now characterize the Atiyah cocycle $\alpha_{\nabla}^E \in \O(A) \otimes A^\perp \otimes \End(E)$ in terms of the brackets $\lambda_\bullet$ coming from $Q_L$. Recall that we started from $D^{L,E} \in (\O(L) \otimes \End(E))^1$ which extends $D^{A,E} \in (\O(A) \otimes \End(E))^1$. This can also be interpreted by a family of degree $1$ linear maps
 $\{\tilde{m}_k:~S^{k-1}(L) \otimes E \rightarrow E\}_{k \geq 1}$ extending $\{m_k:S^{k-1}(A) \otimes E \rightarrow E\}_{k\geq 1}$ (see Equation~\eqref{temp6} for the relation between $D^{L,E}$ and $\tilde{m}_\bullet$).

We further assume that
\bd
\alpha_{\nabla}^E= \sum_{k \geq 0}\alpha_k,
\ed
where $\alpha_k \in S^k(A^\vee) \otimes (L/A)^\vee \otimes \End(E)$. Below is the explicit formula of $\alpha_k$.
\begin{Pro}\label{Atiyah cocycle formula}
   For all $a_1,\cdots, a_k \in A, e \in E, b \in B = L/A$, we have
\begin{align*}
  &\quad(-1)^{k+1}\alpha_k(a_1,\cdots,a_k,b,e) \\
  &= \sum_{p=0}^{k}\sum_{\sigma \in \sh(p,k-p)}\epsilon(\sigma)\tilde{m}_{k-p+3}(\lambda_p(a_{\sigma(1)},\cdots,a_{\sigma(p)}),\cdots, a_{\sigma(k)},l,e) \\
  &\quad+ \sum_{p=0}^{k}\sum_{\sigma \in \sh(p,k-p)}(-1)^{\abs{b}(\ast_k-\dagger_{p}^\sigma)}\epsilon(\sigma)\tilde{m}_{k-p+2}(\lambda_{p+1}(a_{\sigma(1)},\cdots, a_{\sigma(p)}, l),a_{\sigma(p+1)},\cdots,a_{\sigma(k)},e) \\
  &\quad+ \sum_{p=0}^k\sum_{\sigma \in \sh(p,k-p)}(-1)^{\dagger_{p}^\sigma}\epsilon(\sigma)m_{p+1}(a_{\sigma(1)},\cdots,a_{\sigma(p)}, \tilde{m}_{k-p+2}(a_{\sigma(p+1)}, \cdots, a_{\sigma(k)}, l,e)) \\
  &\quad+ \sum_{p=0}^k\sum_{\sigma \in \sh(p,k-p)}(-1)^{\dagger_p^\sigma + \abs{b}(\ast_k - \dagger_p^\sigma+1)}\epsilon(\sigma)\tilde{m}_{p+2}(a_{\sigma(1)}, \cdots, a_{\sigma(p)}, l, m_{k-p+1}(a_{\sigma(p+1)}, \cdots, a_{\sigma(k)},e)),
\end{align*}
where $l\in L$ satisfies $p(l)=b$,  $\ast_k = \sum_{i=1}^k\abs{a_i}$ and $\dagger_p^\sigma = \sum_{i=1}^p\abs{a_{\sigma(i)}}$.
\end{Pro}
The proof follows from some straightforward computations and thus is omitted.
\begin{Rm}\label{Rm:Lsplits}
%Prior to the construction of
To construct the Atiyah cocycle $\alpha_{\nabla}^E$, we need $D^{L,E}$, or $\tilde{m}_k:~S^{k-1}(L)\otimes E\rightarrow E,k\geq 1$. Nevertheless, Proposition \ref{Atiyah cocycle formula} implies that the only information we need  is the behavior of $\tilde{m}_k$ restricted to $S^{k-2}(A)\otimes L\otimes E$. In other words, to compute $\alpha_{\nabla}^E$, it is enough to do first order extensions $\tilde{m}^{(1)}_k:S^{k-2}(A) \otimes L \rightarrow \End(E)$ of $m^E_k$, for all $k \geq 2$. For this reason, we believe that there should exist other, perhaps ``higher'' Atiyah classes.

A more convenient way to get $\alpha_{\nabla}^E$ is to find a complementary subspace to $A$ in $L$. In doing so, one may simply assume that $L=A\oplus B$, where $B$ is only a sub-vector space, not necessarily a subalgebra of $L$. Then $\O(L) \cong \O(A) \otimes \O(B)$. Let $(E,D^{A,E}\sim m^E_\bullet)$ be an $A$-module, where $D^{A,E} \in \O(A) \otimes \End(E) \subset \O(L) \otimes \End(E)$. Thus $\nabla = Q_L + D^{A,E}$ is an $L$-connection on $E$ extending $(A,E)$. Equivalently, $\nabla$ is determined by $\{\tilde{m}^E_k\}_{k\geq 1}: S^{k-1}(L) \otimes E \rightarrow E$:
$$
\tilde{m}^E_{k}= \sum_{p\geq 0}\tilde{m}^E_k\mid_{S^{k-1-p}(A)\otimes S^p(B) \otimes E} =
\begin{cases}
m^E_k, & p=0, \\
0, & p > 0.
\end{cases}
$$
Then the Atiyah cocycle becomes much simpler:
\begin{align*}
  &\quad(-1)^{k+1}\alpha_k(a_1,\cdots,a_k,b,e) \\
  &=\sum_{p=0}^{k}\sum_{\sigma \in \sh(p,k-p)}(-1)^{\abs{b}(\ast_k-\dagger_{p}^\sigma)}\epsilon(\sigma)m^E_{k-p+2}
  (
  	{\prA}(
  			\lambda_{p+1}
  				(
  				a_{\sigma(1)},\cdots, a_{\sigma(p)}, b
  				)
  			),
  			a_{\sigma(p+1)},\cdots,a_{\sigma(k)},e
  ),
\end{align*}
where ${\prA}:L \rightarrow A$ is the projection.
\end{Rm}

From now on, when we talk about the Atiyah cocycle of an SH Lie pair
$(L,A)$ with respect to an $A$-module $E$, we always assume a splitting of sequence~\eqref{SES de GVS} and that the Atiyah cocycle is obtained via the trivial $L$-connection on $E$ extending $(A,E)$ as in Remark~\ref{Rm:Lsplits}.

\begin{Ex}\label{Lie pair}
 Let $(\mathfrak{g},\mathfrak{h})$ be an ordinary Lie algebra pair and $E$ an $\mathfrak{h}$-module, where $\mathfrak{g}, \mathfrak{h}$ and $E$ are all usual ungraded vector spaces. The Atiyah class in \cite{CSX} can be recovered as follows: In fact, setting $L=\mathfrak{g}[1], A=\mathfrak{h}[1]$, we get an  SH Lie pair $(L,A)$ with the obvious $A$-module structure on $E$. Applying Proposition \ref{Atiyah cocycle formula}, we get the Atiyah cocycle
$$
\alpha_{\nabla}^E(a,b,e)= \nabla_{[a,l]}(e) - \nabla_a\nabla_{l}(e) + \nabla_{l}\nabla_a(e)= -R^\nabla(a,l)(e),$$
 where $a \in A, b \in L/A, e \in E$,  $l\in L$ such that $p(l)=b$ and $\nabla:L\otimes E\rightarrow E$ is an $L$-connection  extending $(A,E)$. Comparing with the Atiyah cocycle defined in \cite{CSX}, the only difference is a minus sign.

A nontrivial example of Atiyah classes of this type can be found in \cite{CCT} (see also \cite{CSX}*{Example 22}).
\end{Ex}

%\subsection{Some special situations}

%\subsubsection{One-term SH Lie pairs}
\begin{Ex}
%A one-term $L_\infty[1]$-algebra is an $L_\infty[1]$-algebra $L$ which is concentrated in degree $(-1)$. The only nontrivial structure map is $\lambda_2$. In fact, $L[-1]$ is an ordinary Lie algebra.
Let $(L=L^{-1},A=A^{-1})$ be a one-term SH Lie pair and $E=\oplus_{n\in \mathbb{Z}} E^n$ be an $A$-module, or a Lie algebra representation up to homotopy~\cite{AC} of $A[-1]$ on $E$. Assume that $L=A\oplus B$, where $B$ is also concentrated in degree $(-1)$. If the $A$-module structure of $E$ is given by $m_k:~S^{k-1}A\otimes E\rightarrow E, k \geq 1$, then the Atiyah cocycle $\alpha^E =\sum_{k\geq 0} \alpha_k\in \O(A)\otimes B^\vee\otimes \End(E)$ is given by %(see Remark \ref{Rm:Lsplits}):
\begin{align*}
     (-1)^{k+1} \alpha_k(a_1,\cdots,a_k,b,e)  =
    \sum_{i=1}^{k}(-1)^{k+i} {m}_{k+1}({\prA}\lambda_{2}(a_{i}, b),\cdots,\hat{a_i},\cdots,e)  ,
\end{align*}
for $a_i \in A, e \in E, b \in B = L/A$.

%This Atiyah class is much more complicated than the case of an ordinary Lie algebra pair with respect to an ordinary Lie algebra module (see Example \ref{Lie pair}).
In particular, if the $A$-module structure on $E$ has only two nontrivial actions $m_1:E\rightarrow E$ and $m_2:A\otimes E\rightarrow E$, then the Atiyah cocycle $\alpha^E = \alpha_1\in A^\vee\otimes B^\vee\otimes \End(E)$ reads
\[
\alpha_1(a,b,e)={m_2}({\prA}\lambda_2(a,b),e).
\]
\end{Ex}

%\subsubsection{Two-term SH Lie pairs}
\begin{Ex}
%A two-term $L_\infty[1]$-algebra is a Lie 2-algebra in the sense of \cite{BC}.
Let $(L=L^{-2}\oplus L^{-1},A=A^{-2}\oplus A^{-1})$ be a Lie 2-algebra~\cite{BC} pair with brackets $\lambda_1,\lambda_2,\lambda_3$ and $E=E^{-2}\oplus E^{-1}$ be an $A$-module. Let us fix a splitting $L=A\oplus B$. The Atiyah cocycle $\alpha^E=\alpha_0+\alpha_1+\alpha_2 $ ($\alpha_i\in S^i (A^\vee)\otimes B^\vee\otimes \End(E)$) is given by:
%\[-\alpha_0(b,e)= {m}_2({\prA} \lambda_1(b),e),\]
\begin{align*}
 \alpha_0(b,e)&= -{m}_2({\prA} \lambda_1(b),e) \\
 \alpha_1(a,b,e)&= (-1)^{\abs{b}\abs{a}} {m}_3({\prA} \lambda_1(b),a,e)+ {m}_2({\prA} \lambda_2(a,b),e)\\
%\end{align*}
%and
%\begin{align*}
\alpha_2(a_1,a_2,b,e) &= -{m}_2({\prA} \lambda_3(a_1,a_2,b),e)
 -(-1)^{\abs{b}\abs{a_2}} {m}_3({\prA} \lambda_2(a_1,b),a_2,e) \\
 &\quad- (-1)^{(\abs{b}+\abs{a_2})\abs{a_1}}{m}_3({\prA} \lambda_2(a_2,b),a_1,e),
\end{align*}
for $a_i \in A, e \in E, b \in B = L/A$.
\end{Ex}

%\subsubsection{The Atiyah class of a DGLA pair}

%The notion of a differential graded Lie algebra (DGLA, for short) is an $L_\infty[1]$-algebra with only nontrivial structure maps $\lambda_1,\lambda_2$. %are one and the same (after degree $1$ shifting).
\begin{Ex}
Let $(L,A,\lambda_1,\lambda_2)$ be a DG Lie algebra pair. Suppose that $L=A\oplus B$. Then the associated
$A$-module structure on $B$ consists of two actions:~$m_1^B$ and $m_2^B$. Assume that $E$ is an $A$-module with only two nontrivial actions from $A$:~$m_1 $ and $m_2 $. Then the Atiyah cocycle has two terms
$$
\alpha^E =\alpha_0+\alpha_1\in (B^\vee\otimes \End(E))~\oplus~(A^\vee\otimes B^\vee\otimes \End(E)),
$$
where
\begin{align*}
  -\alpha_0(b,e) &= {m}_2({\prA} \lambda_1(b),e), & \alpha_1(a,b,e) &=  {m}_2({\prA} \lambda_2(a,b),e).
\end{align*}
%Similar to that of Section \ref{Sec:spectral},
For the $A$-module $F=B^\vee\otimes \End(E)$,  we are able to split  the differential operator
$$\partial_A=\partial_0 + \partial_1:~\O(A)\otimes F\rightarrow \O(A)\otimes F,$$
where
\begin{align*}
\partial_0:~&S^{\bullet}(A^\vee)\otimes F\rightarrow S^{\bullet}(A^\vee)\otimes F, &
\partial_1:~&S^{\bullet}(A^\vee)\otimes F\rightarrow S^{\bullet+1}(A^\vee)\otimes F.
\end{align*}
Now the Chevalley-Eilenberg cochain complex $(\O(A)  \otimes F, \partial_A)$ associated to the $A$-module $F$ becomes a double complex:
$$
 D^{p,q}= (S^p(A^\vee)\otimes F)^{p+q},\qquad p\geq 0, q\in \Z
$$
with differentials
$$\partial_0:~D^{p,q}\rightarrow D^{p,q+1},\qquad \partial_1:~D^{p,q}\rightarrow D^{p+1,q}.$$

As for the Atiyah cocycle $\alpha^E=\alpha_0+\alpha_1$,  it sits in $D^{0,2}\oplus D^{1,1}$.  So, the Atiyah class yields two canonical elements:~$[\alpha_0]$ in $\Ha^2(D^{0,\bullet},\partial_0)$ and $[\alpha_1]$ in $\Ha^1(D^{\bullet,1},\partial_1)$.
\end{Ex}

We present a particular example with nontrivial Atiyah classes.
\begin{Ex}
  Let $A = \operatorname{span}\{a_1,a_2\}$ be a $2$-dimensional vector space concentrating in degree $(-1)$ such that $A^\vee = \operatorname{span}\{a_1^\vee,a_2^\vee\}$, and $B = \operatorname{span}\{b\}$ an ordinary $1$-dimensional vector space concentrating in degree $0$
  with dual space $B^\vee = \operatorname{span}\{b^\vee\}$. Then $L = A \oplus B$ together with the homological vector field $Q_L = \delta:~A^\vee \rightarrow S^2(A^\vee) \otimes B^\vee$ defined by
  $$
\delta(a_1^\vee) = k_1 a_1^\vee \odot a_2^\vee \otimes b^\vee, \qquad \delta(a_2^\vee) = k_2 a_1^\vee \odot a_2^\vee \otimes b^\vee, \;\;k_1,k_2 \in \k
  $$
  determine an SH Lie pair $(L,A)$ such that $A \subset L$ is abelian. Note that the only nontrivial structure map is $\lambda_3:~S^2(A) \otimes B \rightarrow A$ by
  \bd
   \lambda_3(a_1,a_2,b) = -k_1a_1 - k_2a_2.
  \ed

  Let $E$ be another $1$-dimensional vector space.
  Then $D^E:~E \rightarrow A^\vee \otimes E$ defined by, for all $e \in E$,
  \bd
    D^E(e) = (k_3a_1^\vee + k_4a_2^\vee) \otimes e, \;\;\text{where}\;k_3,k_4 \in \k\;\text{such that}\;k_1k_3 + k_2k_4 \neq 0,
  \ed
  determines an $A$-module structure on $E$. Equivalently, we have
  $$
 m_2(a_1,e) = -k_3e, \qquad m_2(a_2,e) = -k_4e.
  $$
  The only nontrivial part of the Atiyah cocycle is
  \bd
   \alpha_2(a_1,a_2,b,e) = -{m}_2({\prA} \lambda_3(a_1,a_2,b),e) = -m_2(-k_1a_1-k_2a_2,e) = -(k_1k_3+k_2k_4)e.
  \ed
  Note that the $A$-module structure $\partial_A^{B^\vee \otimes \End(E)}$ on $B^\vee \otimes \End(E)$ is trivial in this case. Thus the Atiyah class $[\alpha_2]\in \Ha_{\mathrm{CE}}^2(A,B^\vee\otimes \End(E))$ is nontrivial.
\end{Ex}
%Another instance with nontrivial Atiyah class is the subsequent Example \ref{Ex:nontrivialalphaB}.

\subsection{An equivalent description of Atiyah classes}\label{Equivalent description}

Let $(L,A)$ be an SH Lie pair. Then there is a coadjoint $A$-module structure $\partial_A^{A^\vee} = Q_A + D^{A^\vee}$ on $A^\vee$, which is the dual of the adjoint $A$-module structure on $A$ (Example~\ref{Adjoint rep.}). Moreover, we have the following commutative diagram:
\[
\xymatrix{
A^\vee \ar[r]^-{Q_A} \ar[d]^-{D^{A^\vee}} & \O(A) \ar[dl]^-{I^A}\\
   \O(A)\otimes A^\vee, &
}
\]
where $I^A$ is defined in Equation~\eqref{I}. In fact, we recall that $Q_A = \sum_{k\geq0}(-1)^{k}\lambda_k^\vee$, and for all $\xi \in A^\vee, a_1\odot\cdots\odot a_k \otimes a_{k+1} \in S^k(A) \otimes A$,
\begin{align*}
  \langle (I^A\circ Q_A)(\xi), a_1\odot\cdots\odot a_k \otimes a_{k+1} \rangle &= (-1)^{k+1}\langle \lambda_{k+1}^\vee(\xi), \mu_{k+1}^A(a_1\odot\cdots\odot a_k\otimes a_{k+1}) \rangle \\
  &= (-1)^{\abs{\xi}+k} \langle \xi, m^A_{k+1}(a_1,\cdots,a_k,a_{k+1}) \rangle,
\end{align*}
where $m^A_{k+1} = \lambda_{k+1} \circ \mu_{k+1}^A: S^k(A) \otimes A \rightarrow A$ is the adjoint $A$-module structure on $A$. Using Equation~\eqref{temp6}, we have
\bd
 \langle (I^A \circ Q_A)(\xi), a_{k+1} \rangle = (-1)^{\abs{\xi}+1} \langle \xi, D^A(a_{k+1}) \rangle = \langle D^{A^\vee}(\xi),a_{k+1} \rangle \in S^k(A^\vee).
\ed
Hence, we have $D^{A^\vee} = I^A \circ Q_A$, as desired.

Similarly, $L$ carries a natural $A$-module structure
$$
m^L_k = \lambda_k \circ \mu^L_{k} \circ (j^{\odot (k-1)}\otimes 1): S^{k-1}A \otimes L \rightarrow L, k \geq 1,
$$
where $\mu^L_{k}$ is the operator defined by Equation~\eqref{mu}. And the dual $A$-module structure $\partial_A^{L^\vee}=Q_A+D^{L^\vee}$ on $L^\vee$ fits into the commutative diagram:
\[
\xymatrix{
L^\vee\ar[r]^-{Q_L} \ar[d]_-{D^{L^\vee}} &\O(L) \ar[dl]^-{J}\\
   \O(A)\otimes L^\vee, &
}
\]
where $J = (j^\vee \otimes 1) \circ I^L$. Moreover, it follows from a simple induction argument that
\begin{align}
\label{kr2} \partial^{L^\vee}_A \circ J &= J \circ Q_L: \O(L) \rightarrow \O(A) \otimes L^\vee, \\
\label{kr3}  \partial^{A^\vee}_A \circ I_A &= I_A \circ Q_A: \O(A) \rightarrow \O(A) \otimes A^\vee.
\end{align}

Using Equations~\eqref{key relation}, \eqref{kr2} and~\eqref{kr3}, it can be verified that the linear dual of Sequence~\eqref{SES de GVS} of graded vector spaces
\bd
\xymatrix{
0 \ar[r] & A^\perp \ar[r] & L^\vee \ar[r]^-{j^\vee} & A^\vee \ar[r] & 0
}
\ed
is also a short exact sequence of $A$-modules.

Let $(E,m_\bullet\sim D^{A,E})$ be an $A$-module. We have a companion exact sequence of $A$-modules:
\bd%\label{SES de modules}
\xymatrix{
0 \ar[r] & A^\perp \otimes \End(E) \ar[r] & L^\vee \otimes \End(E) \ar[r]^-{j^\vee\otimes 1} & A^\vee \otimes \End(E) \ar[r] & 0,
}
\ed
as well as a short exact sequence of cochain complexes:
\be\label{SES}
\xymatrix{
0 \ar[r] & \O(A)\otimes  A^\perp \otimes \End(E)  \ar[r] & \O(A)\otimes L^\vee \otimes \End(E)  \ar[r]^-{1\otimes j^\vee\otimes 1} & \O(A)\otimes A^\vee \otimes \End(E) \ar[r] & 0.
}
\ee
A long exact sequence on the cohomology level follows:
\begin{eqnarray}\nonumber
\cdots &\rightarrow \Ha^1_{\mathrm{CE}}(A,A^\perp \otimes \End(E))\longrightarrow \Ha^1_{\mathrm{CE}}(A,L^\vee\otimes \End(E))\longrightarrow
\Ha^1_{\mathrm{CE}}(A,A^\vee\otimes \End(E))\\
\label{LES}
&\xrightarrow{\delta}\Ha^2_{\mathrm{CE}}(A,A^\perp\otimes \End(E))\longrightarrow \Ha^2_{\mathrm{CE}}(A,L^\vee\otimes \End(E))\longrightarrow\cdots.
\end{eqnarray}

\begin{Lem}
The element  $(I^A \otimes 1)(D^{A,E})\in \O(A)\otimes A^\vee \otimes \End(E)$ is a degree $1$ cocycle.
\end{Lem}
\bp
Since $I^A$ is a derivation valued in the $\O(A)$-bimodule $\O(A) \otimes A^\vee$, $(I^A \otimes 1)$ is also a derivation on the associative algebra $(\O(A) \otimes \End(E),\circ)$. Thus
\[
(I^A\otimes 1)((D^{A,E})^2)=(I^A\otimes 1)(D^{A,E})\circ D^{A,E} + D^{A,E} \circ (I^A\otimes 1)(D^{A,E}) = [D^{A,E},(I^A\otimes 1)D^{A,E}].
\]

Using Equation~\eqref{kr3}, we have
\bd
 (\partial_A^{A^\vee} \otimes 1)((I^A\otimes 1)(D^{A,E}))=(I^A\otimes 1)(Q_A(D^{A,E})).
\ed
Hence,
\begin{align*}
&\quad \partial_A^{A^\vee \otimes \End(E)}((I^A \otimes 1)(D^{A,E}))=(\partial_A^{A^\vee} + [D^{A,E},-])((I^A\otimes 1)(D^{A,E})) \\
 &=(I^A\otimes 1)(Q_A(D^{A,E})) + (I^A\otimes 1)((D^{A,E})^2) = (I^A \otimes 1)(Q_A(D^{A,E})+(D^{A,E})^2)=0,
\end{align*}
where the last equality follows from the Maurer-Cartan equation~\eqref{MC-EQ}.
\ep
It turns out that the element $ (I^A \otimes 1)(D^{A,E})$ gives the Atiyah class:
\begin{Thm}\label{Thm:deltagivesatyah}
The cohomology class
\[
\delta[(I^A \otimes 1)(D^{A,E})] \in \Ha^2_{\mathrm{CE}}(A,A^\perp \otimes \End(E))
\]
coincides with the Atiyah class $[\alpha^E]$.
\end{Thm}
\bp
We chase the connecting map $\delta$ in Equation~\eqref{SES}: Starting with $(I^A \otimes 1)(D^{A,E})$, one first chooses a degree $1$ element
$\beta\in \O(A)\otimes L^\vee\otimes \End(E)$ such that $(1\otimes j^\vee \otimes 1)(\beta)=(I^A \otimes 1)(D^{A,E}).$

Then one is able to find a unique degree $2$ element $\alpha \in \O(A)\otimes A^\perp \otimes \End(E)$ such that
\[
\partial_A^{L^\vee\otimes \End(E)} (\beta)=(Q_A+ D^{L^\vee} + [D^{A,E},-])(\beta) = \alpha.
\]
The cohomology class $[\alpha]$ is the upshot of $\delta[(I^A \otimes 1)(D^{A,E})]$. We now show that there exists an $L$-connection $\nabla=Q_L+D^{L,E}$ extending $(A,E)$, i.e., $(j^\vee\otimes 1)(D^{L,E})=D^{A,E}$, and the resulting Atiyah cocycle $\alpha_{\nabla}^E$ equals $\alpha$.

First of all, there exists an element $D^{L,E} \in \O(L) \otimes \End(E)$ such that
\begin{align*}
  (J \otimes 1)(D^{L,E}) = \beta, & \qquad (j^\vee \otimes 1)(D^{L,E}) = D^{A,E}.
\end{align*}
In fact, as $J:\O(L) \rightarrow \O(A)\otimes L^\vee$ is surjective, we can find some $K^{L,E}\in \O(L) \otimes \End(E)$ such that
$(J\otimes 1)(K^{L,E})=\beta$. Then by Equation~\eqref{IJ}, we have
$$
(I^A \otimes 1)(j^\vee\otimes 1)(K^{L,E})=(1\otimes j^\vee\otimes 1)(J\otimes 1)(K^{L,E}) =(1\otimes j^\vee\otimes 1)(\beta)=(I^A \otimes 1)(D^{A,E}).
$$
Note that $\ker(I^A) \cong \ker(r^+)=\k$. Thus $(j^\vee\otimes 1)(K^{L,E}) - D^{A,E} = \varphi$ for some $\varphi \in \End(E)$. It follows that $D^{L,E} = K^{L,E} - \varphi$ satisfies the above requirements.

Then, using Equation~\eqref{kr2},
\begin{align*}
  \alpha_{\nabla}^E &=(J \otimes 1)(Q_L(D^{L,E})+(D^{L,E})^2) = (\partial_A^{L^\vee} \otimes 1) \circ (J \otimes 1)(D^{L,E}) + [(j^\vee \otimes 1)D^{A,E}, (J\otimes 1)(D^{L,E})]  \\
  &=(Q_A + D^{L^\vee})(J\otimes 1)(D^{L,E})+[D^{A,E},(J\otimes 1)(D^{L,E})]\\
  &=(Q_A + D^{L^\vee} + [D^{A,E},-])(\beta) = \alpha,
\end{align*}
as required.
\ep

\subsection{Vanishing of Atiyah classes}

Let $A$ be an $L_\infty[1]$-algebra and $B$ an $A$-module. Then the associated abelian extension $L = A \oplus B$ of $A$ along $B$ (see Proposition~\ref{Thm:module}) gives rise to an SH Lie pair $(L,A)$, while the Atiyah class $\alpha_E$ with respect to any $A$-module $E$ is always trivial. So apparently the Atiyah class measures the nontriviality of the extension of $A$ to $L$. %However, the converse statement is NOT true.

It is natural to ask what we can say in general when the Atiyah class vanishes. The following facts are some first stage results. Further investigations of this question will be shown somewhere else.

\begin{Thm}\label{Vanishing theorem}
  Let $(L,A)$ be an SH Lie pair and $(E,(\partial_A^E=Q_A + D^{A,E}) \sim f_\bullet^E)$ an $A$-module. Then the following four statements are equivalent:
  \begin{enumerate}
    \item The Atiyah class $[\alpha^E] \in \Ha^2_{\mathrm{CE}}(A,A^\perp \otimes \End(E))$ vanishes;
    \item There exists a degree $1$ cocycle $\phi \in \O(A)\otimes L^\vee\otimes \End(E)$ such that
      \be\label{conditiononphi}
        (1\otimes j^\vee\otimes 1)(\phi)=(I^A \otimes 1)(D^{A,E});
      \ee
    \item There exists an $A$-module morphism $\{\phi_k: S^k(A) \otimes L \rightarrow \End(E)[1]\}_{k\geq0}$ from $L$ to $\End(E)[1]$ extending the canonical $A$-module morphism $\phi^{A,E}$ defined in Equation~~\eqref{moduleasmorphism} from $A$ to $\End(E)[1]$, i.e.,
     \be\label{phim}
        \phi_k \circ (1 \otimes j) = \phi_{k}^{A,E} = f^E_{k+1} \circ \mu_{k+1}^A: S^k(A) \otimes A \rightarrow \End(E)[1];
     \ee
%     where $\mu_{k+1}^A: S^k(A) \otimes A \rightarrow S^{k+1}(A)$ is the map defined by~\eqref{mu}.
    \item There exists an $L$-connection $\nabla$ on $E$ extending $(A,E)$ such that the Atiyah cocycle $\alpha_\nabla^E$ of $E$ relative to $\nabla$ vanishes.
  \end{enumerate}
\end{Thm}
\bp
$(1) \Rightarrow (2)$. Assume that $[\alpha^E] = 0$. It follows from Theorem~\ref{Thm:deltagivesatyah} that $\delta[(I^A \otimes 1)(D^{A,E})] = 0$. By chasing the long exact sequence~\eqref{LES}, there exists some $\partial_A^{L^\vee \otimes \End(E)}$-cocycle $\tilde{\phi} \in \O(A)\otimes L^\vee\otimes \End(E)$ of degree $1$ such that
\bd
   [(1\otimes j^\vee\otimes 1)(\tilde{\phi})] = [(I^A \otimes 1)(D^{A,E})] \in \Ha^1_{\mathrm{CE}}(A,A^\vee \otimes \End(E)).
\ed
It follows that there is a degree $0$ element $\beta \in \O(A) \otimes A^\vee \otimes \End(E)$ such that
\bd
 (1\otimes j^\vee\otimes 1)(\tilde{\phi}) - (I^A \otimes 1)(D^{A,E}) = \partial_A^{A^\vee \otimes \End(E)}(\beta).
\ed
By exactness of Sequence~\eqref{SES}, one can choose an element $\gamma \in (\O(A) \otimes L^\vee \otimes \End(E))^0$ such that $(1\otimes j^\vee\otimes 1)(\gamma) = \beta$.

Let $\phi = \tilde{\phi} - \partial_A^{L^\vee \otimes \End(E)}(\gamma) \in (\O(A)\otimes L^\vee\otimes \End(E))^1$. Then
\begin{align*}
 %= (1\otimes j^\vee\otimes 1)(\tilde{\phi} - \partial_A^{L^\vee \otimes \End(E)}(\gamma))
 (1\otimes j^\vee\otimes 1)(\phi) &= (1\otimes j^\vee\otimes 1)(\tilde{\phi})-((1\otimes j^\vee\otimes 1)\circ\partial_A^{L^\vee \otimes \End(E)})(\gamma) \\
 &= (I^A \otimes 1)(D^{A,E}) + \partial_A^{A^\vee \otimes \End(E)}(\beta) - (\partial_A^{A^\vee \otimes \End(E)}\circ (1\otimes j^\vee\otimes 1))(\gamma) \\
 &= (I^A \otimes 1)(D^{A,E}).
\end{align*}

$(2) \Leftrightarrow (3)$. Note that a degree $1$ element $\phi \in \O(A) \otimes L^\vee \otimes\End(E)$ consists of a family of degree $0$ maps $\phi_k: S^k(A) \otimes L \rightarrow \End(E)[1]$. It can be easily seen that Equation~\eqref{conditiononphi} is equivalent to Equation~\eqref{phim}.

$(2) \Rightarrow (4)$. Given a degree $1$ cocycle $\phi \in \O(A) \otimes L^\vee \otimes \End(E)[1]$ satisfying Equation~\eqref{conditiononphi}, by the argument in the proof of Theorem~\ref{Thm:deltagivesatyah}, we can find $D^{L,E} \in \O(L) \otimes \End(E)$ such that
\begin{align*}
  (J \otimes 1)(D^{L,E}) = \phi, & \qquad (j^\vee \otimes 1)(D^{L,E}) = D^{A,E}.
\end{align*}
Then $\nabla = Q_L + D^{L,E}$ is an $L$-connection on $E$ extending $(A,E)$. The associated Atiyah cocycle
\begin{align*}
  \alpha_{\nabla}^E &= (J\otimes 1)(R^\nabla) = (J \otimes 1)(Q_L(D^{L,E}) + D^{L,E} \circ D^{L,E}) \\
                   &= \partial_A^{L^\vee}((J \otimes 1)(D^{L,E})) + [(j^\vee\otimes 1)(D^{L,E}), (J\otimes 1)(D^{L,E})] \qquad (\text{by Equations~\eqref{kr2} and~\eqref{J1}})\\
                   &= \partial_A^{L^\vee}(\phi) + [D^{A,E}, \phi] = \partial_A^{L^\vee \otimes \End(E)[1]}(\phi) = 0.
\end{align*}
Finally, $(4) \Rightarrow (1)$ is obvious. This completes the proof.
\ep

\section{Atiyah classes as functors}\label{Atiyah as functors}
\subsection{Atiyah operators}
%We maintain all previous settings.
Let $(L,A)$ be an SH Lie pair and denote by $B=L/A$ the standard $A$-module. The identification of $B^\vee$ and $A^\perp$ is assumed. For simplicity, from now on, we will denote the Chevalley-Eilenberg differential $\partial_A^E$ of any $A$-module $E$ by $\partial_A$.

Recall that the Atiyah cocycle
$$
\alpha^E_\nabla\in \O(A)\otimes A^\perp \otimes \End(E)=\O(A)\otimes B^\vee\otimes \Hom(E,E)
=\O(A)\otimes \Hom(B\otimes E,E),
$$
where $\nabla$ is an $L$-connection on $E$ extending $(A,E)$. The associated $\O(A)$-linear maps
\begin{align*}
\bm{\alpha}_\nabla^E:\qquad &\O(A)\otimes E\longrightarrow \O(A)\otimes B^\vee\otimes E;\\
\bm{\alpha}_\nabla^E(x):\qquad &\O(A)\otimes E\rightarrow \O(A)\otimes E,\qquad \mbox{ where } x\in \O(A)\otimes B,
%\overrightarrow{\alpha_\nabla^E}:\qquad & \O(A)\otimes B\longrightarrow \Hom_{\O(A)}(\O(A)\otimes E,\O(A)\otimes E)\cong \O(A)\otimes \End(E),
\end{align*}
will be called Atiyah operators. As $\alpha^E_\nabla$ is a cocycle, we have
\begin{Lem}
If $x\in \O(A)\otimes B$ is a cocycle, i.e. $\partial_A(x)=0$, then $\partial_A^{\End(E)}(\bm{\alpha}^E_{\nabla}(x)) = 0$, %the Atiyah operator $\bm{\alpha}^E_{\nabla}(x)$ is a chain map,
i.e., $\bm{\alpha}^E_{\nabla}(x)$ is an $A$-module morphism from $E$ to itself.
\end{Lem}
%As in Remark \ref{Rm:Lsplits}, a splitting of the vector space $L\cong A\oplus B$ would simplify a lot in computing Atiyah cocycles. So we directly assume that $L=A \oplus B$, though it is not canonical. We then have
Let us fix a splitting of the short exact sequence~\eqref{SES de GVS} so that $L \cong A \oplus B$ and
\[\O(L)\cong \O(A)\odot \O(B)\cong \O(A)\otimes \O(B).\]
And the associated homological vector field $Q_L\in \Der(\O(L))$ decomposes into a sum of derivations
\be\label{T}
Q_L=Q_A+\delta+\mathcal{R}+D^\perp+\sum_{i\geq 2}T_i,
\ee
%Here $Q_A$ is the homological vector field of $A$, which is determined by the restriction of $Q_L$ on $A^\vee$ and then  projected to  $\O(A)$. The derivation $\delta$ is determined by the restriction of $Q_L$ on $A^\vee$ and then  projected to  the $\O(A)\otimes B^\vee$-part. The derivation $\mathcal{R}$ is determined by the restriction of $Q_L$ on $A^\vee$ and then  projected to  the $\O(A)\otimes  \widehat{S}{^{\geq 2}(B^\vee)}$-part. When $Q_A$, $\delta$ and $\mathcal{R}$ act on $B^\vee$, they only yield zeros.  The derivation $D^\perp$ is the $A$-module structure on $B^\vee$. In other words, $D^\perp$ is determined by the restriction of $Q_L$ on $B^\vee$ and then  projected to  the $\O(A)\otimes B^\vee$-part. Similarly, $T_i$ is determined by the restriction of $Q_L$ on $B^\vee$ and then  projected to  the $\O(A)\otimes S^i(B^\vee)$-part. Also, $D^\perp$ and $T_i$ act trivially on $A^\vee$.

%We summarize these operators in the following table:
where
$$
\begin{cases} Q_A &= Q_L|_{:A^\vee\longrightarrow \O(A)};\\
\delta &= Q_L|_{:A^\vee\longrightarrow \O(A)\otimes B^\vee};\\
\mathcal{R} &= Q_L|_{:A^\vee\longrightarrow \O(A)\otimes  \widehat{S}{^{\geq 2} (B^\vee)}};\\
D^\perp &= Q_L|_{:B^\vee\longrightarrow \O(A)\otimes B^\vee};\\
T_i &= Q_L|_{:B^\vee\longrightarrow \O(A)\otimes S^i(B^\vee),}\quad i\geq 2,
\end{cases}
$$
and they are extended as a derivation of $\O(L)$ in a natural manner.

In this situation, there is an $L$-connection $\nabla=Q_L+D^{A,E}$ extending $(A,E)$ (see Remark~\ref{Rm:Lsplits}). The associated Atiyah cocycle is denoted by $\alpha^E$, where the subscript is omitted. Similarly, the associated Atiyah operator will be denoted by $\bm{\alpha}^E$.

We extend the operator $\delta$ to a $\k$-linear and degree $1$ map for any graded vector space $E$
\begin{align*}
\delta=(\delta\otimes 1):\O(A)\otimes E&\rightarrow \O(A)\otimes B^\vee \otimes E,%\\  \omega\otimes e&\mapsto \delta(\omega)\otimes e,
\end{align*}
such that the Leibniz rule
\bd
 \delta(\xi \odot \eta \otimes e) = (-1)^{(\abs{\xi}+1)\abs{\eta}}\eta \odot \delta(\xi) \otimes e + (-1)^{\abs{\xi}}\xi \odot \delta(\eta) \otimes e
\ed
holds for all $\xi,\eta \in \O(A), e \in E$.

\begin{Lem}\label{Important}
As a map $\O(A)\otimes E\rightarrow \O(A)\otimes B^\vee \otimes E$, the Atiyah operator
$$
\bm{\alpha}^E=[\delta,\partial_A]=\delta\circ \partial_A+\partial_A\circ \delta.
$$
\end{Lem}
\bp
%By the decomposition \eqref{T}, when $Q_L$ is restricted to   $\O(A)\otimes B^\vee$ and then projected to  $\O(A)\otimes B^\vee$, it is simply the sum of
%  $ Q_A$ and $D^\perp$.
Recall the definition of Atiyah cocycles:
\begin{align*}
\alpha^E &= (J\otimes 1)(\nabla^2)= (J\otimes 1)\left(Q_L(D^{A,E}) + (D^{A,E})^2\right) \\
         &= \partial_A^{L^\vee}((J \otimes 1)(D^{A,E})) + [D^{A,E},(J\otimes1)(D^{A,E})]\;\qquad\qquad\quad (\text{by Equations~\eqref{kr2} and~\eqref{J1}}) \\
         &= \delta(D^{A,E}) + \partial^{A^\vee}_A((I^A \otimes 1)(D^{A,E}))+ [D^{A,E},(I^A\otimes1)(D^{A,E})] \quad\qquad\quad(\text{by Equation~\eqref{MC-EQ}}) \\
 &=\delta(D^{A,E}) + (I^A \otimes 1)(Q_A(D^{A,E}) + \frac{1}{2}[D^{A,E}, D^{A,E}])\qquad\quad(\text{by Equations~\eqref{kr3} and~\eqref{derivation of I}}) \\
         &= \delta(D^{A,E}) \in \O(A) \otimes B^\vee \otimes \End(E).
\end{align*}
By observing $Q_L^2 = 0$ on the $\O(A)$ to $\O(A) \otimes B^\vee$-part, we have
\bd
 \partial_A^\perp \circ \delta + \delta \circ Q_A = 0: \O(A) \rightarrow \O(A) \otimes B^\vee.
\ed
Thus, for all $\omega \otimes e \in \O(A) \otimes E$,
\begin{align*}
  &\quad\partial_A(\delta(\omega \otimes e)) + \delta(\partial_A(\omega \otimes e)) \\
  &= \partial_A^\perp(\delta(\omega)) \otimes e + (-1)^{\abs{\omega}+1}\delta(\omega)\odot(D^{A,E}(e)) + \delta(Q_A(\omega))\otimes e + (-1)^{\abs{\omega}}\delta(\omega \odot D^{A,E}(e)) \\
  &= \omega \odot \delta(D^{A,E})(e) = \delta(D^{A,E})(\omega \otimes e) = \bm{\alpha}^E(\omega \otimes e).
\end{align*}

\ep
Applying Lemma \ref{Important}, one can easily get the following properties of Atiyah operators:
\begin{Lem}\label{properties of Atiyah operators}
  For all $x \in \O(A)\otimes B$, the Atiyah operator $\bm{\alpha}^\bullet(x)$ has the following properties:
  \begin{enumerate}
    \item For any $A$-modules $E$ and $F$,
         $$
        \bm{\alpha}^{E\otimes F}(x)(r\otimes_{\O(A)} s)=(\bm{\alpha}^E(x)r)\otimes_{\O(A)} s+(-1)^{\abs{x}\abs{r}}r\otimes_{\O(A)} (\bm{\alpha}^F(x)s),
         $$
         for all $r\in \O(A)\otimes E$, $s\in \O(A)\otimes F$;
    \item For any $A$-module $E$ with dual module $E^\vee$, and for all $\varphi\in \O(A)\otimes E^\vee, r\in \O(A)\otimes E$,
        $$
         \langle \bm{\alpha}^{E^\vee}(x)\varphi,r\rangle = -(-1)^{\abs{x}\abs{\varphi}}\langle \varphi, \bm{\alpha}^{E}(x)r\rangle; %\qquad \forall \varphi\in \O(A)\otimes E^\vee, r\in \O(A)\otimes E.
        $$
    \item For $A$-modules $E$ and $F$,
        $$
        (\bm{\alpha}^{\Hom(E,F)}(x)\kappa)(r)=\bm{\alpha}^F(x)(\kappa(r))-(-1)^{\abs{x}\abs{\kappa}} \kappa(\bm{\alpha}^E(x)r),
        $$
        for all $\kappa\in \O(A)\otimes \Hom(E,F)$, $r\in \O(A)\otimes E$.
  \end{enumerate}
\end{Lem}

\subsection{Atiyah classes as Lie structures}
Let $(L,A)$ be an SH Lie pair, and suppose that $L=A\oplus B$ as graded vector spaces, and $E$ an $A$-module. Note that the Atiyah cocycle $\alpha^E$ is a degree 2 element in
$\O(A)\otimes B^\vee \otimes \End(E)$. If we set $\mathbb{B}=B[-2]$, then the associated Atiyah operators are of degree $0$:
\begin{align*}
\bm{\alpha}^B(-)-:\qquad &(\O(A)\otimes \mathbb{B}) \otimes_{\O(A)} (\O(A)\otimes \mathbb{B}) \rightarrow \O(A)\otimes \mathbb{B};\\
\bm{\alpha}^E(-)-:\qquad &(\O(A)\otimes \mathbb{B}) \otimes_{\O(A)} (\O(A)\otimes E) \rightarrow \O(A)\otimes E.
\end{align*}

Here are the main results in this section:
%Since $\alpha^B$ and $\alpha^E$ are both cocycles, the corresponding $\bm{\alpha}^B$ and $\bm{\alpha}^E$ are all chain maps.
\begin{Thm}\label{Cor:LieandLieModuleStr}
Let $(L,A)$ be an  SH Lie pair with the quotient space $L/A=B$. Then the graded vector space $\Ha^\bullet_{\mathrm{CE}}(A,\mathbb{B})$ with the binary operation induced by the Atiyah operator $\bm{\alpha}^B$ is a  Lie algebra. Furthermore, if $E$ is an $A$-module, then $\Ha^\bullet_{\mathrm{CE}}(A,E)$  is a Lie algebra module over $\Ha_{\mathrm{CE}}^\bullet(A,\mathbb{B})$, with the action induced by the Atiyah operator $\bm{\alpha}^E$.

In particular, $\Ha^0_{\mathrm{CE}}(A,\mathbb{B})$  is an ordinary Lie algebra and $\Ha^0_{\mathrm{CE}}(A,E)$ is an ordinary Lie algebra module over $\Ha_{\mathrm{CE}}^0(A,\mathbb{B})$.
\end{Thm}
\begin{Rm}
By Theorem-Definition~\ref{main prop-def}, the Lie algebra and Lie algebra module structures on the cohomology level are all canonical, i.e., they do not depend on the choice of the splitting $L=A\oplus B$.
\end{Rm}

%Below are the proof of Theorem~\ref{Cor:LieandLieModuleStr}:
To proceed the proof, we need some preparations.
Let %$\tau$ be the symmetric exchanging operator:
$$
\tau:B^\vee\otimes B^\vee \rightarrow B^\vee\otimes B^\vee, \qquad \xi\otimes \eta\mapsto (-1)^{\abs{\xi}\abs{\eta}}\eta\otimes \xi.
$$
\begin{Lem}\label{Lem:PreSkewSymm}
The symmetrization of the Atiyah cocycle $\alpha^{B^\vee}$  vanishes up to homotopy, i.e.,
\begin{equation}\label{Eqt:skewsym}
(1\otimes \tau)\alpha^{B^\vee}+\alpha^{B^\vee}=\partial_A P,
\end{equation}
for some $P\in \O(A)\otimes B^\vee\otimes B^\vee\otimes B$.
\end{Lem}
\bp
First, observe the following commutative diagram
\[
\xymatrix{
    B^{\vee}\otimes B^\vee \ar[rr]^-{\frac{1}{2}(1+\tau)} \ar[drr]^-{s} & & B^{\vee}\otimes B^\vee \ar@<.7ex>[d]^-{s} \\
    && S^2(B^{\vee}) \ar@<.7ex>[u]^{s^{-1}},
}
\]
where the operations $s$ and $s^{-1}$ are defined by
$$
s:~\xi\otimes \eta\mapsto \xi\odot \eta,\qquad s^{-1}:\xi\odot \eta\mapsto \frac{1}{2}(1+\tau)(\xi\otimes \eta),\quad \forall \xi,\eta\in B^\vee.
$$
It is clear that $s^{-1}$ is right inverse of the symmetrization operator $s$, i.e., $s \circ s^{-1} = id:~S^2(B^\vee) \rightarrow S^2(B^\vee)$. We introduce
$$
\delta^0:\O(A)\otimes B^\vee\rightarrow \O(A)\otimes S^2(B^\vee),\qquad \delta^0=(1\otimes s)\circ \delta.
$$
Then we get the following commutative diagram
\[
\xymatrix{
   B^{\vee} \ar[rrr]^-{\frac{1}{2}(\bm{\alpha}^{B^\vee}+(1\otimes \tau)\bm{\alpha}^{B^\vee})} \ar[drrr]^-{\delta^0\circ D^\perp} & & & \O(A) \otimes (B^{\vee}\otimes B^\vee) \ar@<3.0ex>[d]^-{1 \otimes s} \\
    &&& \O(A) \otimes S^2(B^{\vee}) \ar@<-2.0ex>[u]^{1 \otimes s^{-1}}.
}
\]
Namely,
$$
(1\otimes s)\circ \frac{1}{2}(\bm{\alpha}^{B^\vee}+(1\otimes \tau)\bm{\alpha}^{B^\vee})=\delta^0\circ D^\perp.
$$
In fact, it amounts to check
$$
\bm{\alpha}^{B^\vee}=\delta\circ D^\perp,\;\text{as a map}\; B^\vee\rightarrow \O(A)\otimes B^\vee\otimes B^\vee,
$$
which follows from Lemma \ref{Important} and $\partial_A^{B^\vee}|_{B^\vee}=D^\perp$. Hence
$$
(1\otimes \tau)\alpha^{B^\vee}+\alpha^{B^\vee}= 2 (1\otimes s^{-1} \otimes 1) (\delta^0\circ D^\perp).
$$
%As $1\otimes s^{-1} \otimes 1$ is a chain map,
To prove Equation~\eqref{Eqt:skewsym}, it suffices to show that
$$
\delta^0\circ D^\perp \in \Hom(B^\vee, \O(A) \otimes S^2(B^\vee)) \cong \O(A)\otimes S^2(B^\vee)\otimes B
$$
is a coboundary.

Restricting the condition $Q_L^2=0$ on the $B^\vee$ to $\O(A) \otimes S^2(B^\vee)$-part, we get
$$
\delta^0\circ D^\perp+Q_A\circ T_2+D^\perp\circ T_2+T_2\circ D^\perp=0,
$$
where $T_2\in \O(A)\otimes S^2(B^\vee)\otimes B$ is defined in Equation~\eqref{T}. Hence, we have
\begin{align*}
\delta^0\circ D^\perp&=-(Q_A+D^\perp)\circ T_2-T_2\circ D^\perp =-[Q_A+D^\perp,T_2]=-\partial_A T_2,
\end{align*}
as desired.
\ep
\begin{Lem}\label{Lem:preAction}
For all $x,y\in \O(A)\otimes B,r\in \O(A)\otimes E$, we have
$$
\bm{\alpha}^E(x)(\bm{\alpha}^E(y)r)-(-1)^{\abs{x}\abs{y}}\bm{\alpha}^E(y)(\bm{\alpha}^E(x)r)=\bm{\alpha}^E(\bm{\alpha}^B(x)y)r+(\partial_A T)\llcorner (x\otimes_{\O(A)} y\otimes_{\O(A)} r),
$$
where $T=\delta(\alpha^E)\in \O(A)\otimes B^\vee\otimes B^\vee\otimes \End(E)$.
\end{Lem}
\bp
Since $\partial_A(\alpha^E)=0$, it follows from~Lemma \ref{Important} that
$$
\partial_A T=\partial_A(\delta(\alpha^E))=(\partial_A \circ \delta + \delta\circ \partial_A)(\alpha^E) = \bm{\alpha}^{B^\vee \otimes \End(E)}(\alpha^E) = \bm{\alpha}^{\Hom(B\otimes E,E)}(\alpha^E).
$$
Applying Lemma~\ref{properties of Atiyah operators}, we have
\begin{align*}
&\qquad \bm{\alpha}^{\Hom(B\otimes E,E)}(\alpha^E) \llcorner (x\otimes_{\O(A)} y\otimes_{\O(A)} r) = \big(\bm{\alpha}^{\Hom(B\otimes E,E)}(x)(\alpha^E)\big)\llcorner (y\otimes_{\O(A)} r)\\
&=\bm{\alpha}^E(x)(\bm{\alpha}^E(y\otimes_{\O(A)} r))-(-1)^{2\abs{x}}\bm{\alpha}^E(\bm{\alpha}^{B\otimes E}(x)(y\otimes_{\O(A)} r))\\ &=\bm{\alpha}^E(x)(\bm{\alpha}^E(y) r)-\bm{\alpha}^E((\bm{\alpha}^B(x)y)\otimes_{\O(A)} r+(-1)^{\abs{x}\abs{y}}y\otimes_{\O(A)} (\bm{\alpha}^E(x)r))\\ &=\bm{\alpha}^E(x)(\bm{\alpha}^E(y) r)-\bm{\alpha}^E(\bm{\alpha}^B(x)y)r-(-1)^{\abs{x}\abs{y}}\bm{\alpha}^E(y)(\bm{\alpha}^E(x)r).
\end{align*}
\ep

We are now ready to turn to
\bp[Proof of Theorem~\ref{Cor:LieandLieModuleStr}]
Lemma~\ref{Lem:PreSkewSymm} implies that the bracket $[x,y]= \bm{\alpha}^B(x)y$ is skew-symmetric on the cohomology level. When $E$ is taken as $B$ in Lemma~\ref{Lem:preAction}, we see that the $[-,-]$-bracket satisfies Jacobi identity on the cohomology level. Again by Lemma \ref{Lem:preAction}, the operation $x\triangleright r=\bm{\alpha}^E(x)r$ defines a Lie algebra action of $\Ha^\bullet_{\CE}(A,\mathbb{B})$ on $\Ha_{\CE}^\bullet(A,E)$.
\ep

We remark that similar results appeared in \cite{CSX} for that of Lie pairs, and in~\cite{Bottacin} of relative Lie algebroids. In the meantime, the work   related to some facts in the derived categories claimed in \cite{CSX14} is still going on.

\begin{Ex}\label{Ex:nontrivialalphaB}
Take $L=A\oplus B$, where $A$ is spanned by three vectors $a_0,a_1$ and $c$, $B$ by one vector $b$. The degrees are assigned:
$$
  \abs{a_0} =\abs{a_1} =\abs{b}=-1, \qquad \abs{c}=0.
$$
Let the dual vectors be $a_i^\vee$, $c^\vee$ and $b^\vee$, with degrees
$$
\abs{a_0^\vee}=\abs{a_1^\vee} =\abs{b^\vee}= 1, \qquad \abs{c^\vee} = 0.
$$
%\begin{align*}
%  \abs{a_0^\vee}&=\abs{a_1^\vee} =\abs{b^\vee}= 1, & \abs{c^\vee} &= 0.
%\end{align*}

The $Q$-structure on $L$ is the sum of three parts
$$Q_L=Q_A+D^\perp+\delta.$$
Here $Q_A$ is determined by
$$
 Q_A(a_0^\vee)= Q_A(c^\vee)=0, \qquad Q_A(a_1^\vee)= -a_0^\vee\odot a_1^\vee.
$$

The $A$-module structure on $B^\vee$ is given by
\begin{align*}
  D^\perp(b^\vee) &= a_0^\vee\otimes b^\vee.% & \partial_A^B(b) &= a_0^\vee \otimes b.
\end{align*}
The $\delta$-operator is given by
$$\delta(\xi)=(-1)^{\abs{\xi}}\Delta(\xi)\otimes b^\vee,\qquad \forall \xi\in \O(A),$$
where $\Delta$ is a degree $0$ derivation on $\O(A)$ determined by
$$
\Delta(a_0^\vee)=a_1^\vee, \qquad \Delta(a_1^\vee) =a_0^\vee\odot c^\vee, \qquad \Delta(c^\vee)=0.
$$

It follow from some direct computations that $(L,A)$ is an SH Lie pair. The Atiyah operator $\bm{\alpha}^{B^\vee}$ now reads
$$\bm{\alpha}^{B^\vee}(b^\vee)=\delta\circ D^\perp(b^\vee)=\delta(a_0^\vee\otimes b^\vee)=-\Delta(a_0^\vee)\otimes b^\vee\otimes b^\vee=-a_1^\vee\otimes b^\vee\otimes b^\vee.$$
Or, the Atiyah cocycle is spelled as $$\alpha^B=-a_1^\vee\otimes b^\vee\otimes b^\vee\otimes b.$$
The Atiyah class $[\alpha^B]\neq 0$. In fact, any attempt to make $[\alpha^B]= 0$
yields the equation
$$Q_A(\xi)+\xi\odot a_0^\vee=-a_1^\vee,$$
for $\xi\in \O(A)$ with $\abs{\xi}=0$. It obviously has no solution.

The space $\O(A)\otimes \mathbb{B}$ is generated by one element $b[-2]$, and the Lie bracket on $\Ha^\bullet(\O(A)\otimes \mathbb{B})$ can be explicitly expressed:
$$
\bigl[b[-2],b[-2]\bigr]=\bm{\alpha}^{B}(b[-2])b[-2]=a_1^\vee\otimes b[-2].
$$
\end{Ex}

\subsection{Atiyah functors}
Let $A$ be an $L_\infty[1]$-algebra. Then taking the Chevalley-Eilenberg cohomology $\Ha^\bullet_{\mathrm{CE}}(A,-)$ defines a functor% from the category $\Mod_A$ of $A$-modules to the category $\text{GVS}_\k$ of $\k$-vector spaces
\bd
 \Ha^\bullet_{\mathrm{CE}}(A;-):~\Mod_A  \rightarrow \text{GVS}_\k,\;  E  \mapsto \Ha^\bullet_{\mathrm{CE}}(A,E) = \Ha^\bullet(\O(A) \otimes E, \partial_A^E).
\ed
% Meanwhile, $\Ha^\bullet_{\mathrm{CE}}(A,E)$ is also a graded module over the graded algebra $\mathcal{H}= \Ha^\bullet_{\mathrm{CE}}(A,\k)$.
For a morphism $\phi \in \Hom_A(E,F)$, the functor sends $\phi$ to $[\phi]: \Ha^\bullet_{\mathrm{CE}}(A,E) \rightarrow \Ha^\bullet_{\mathrm{CE}}(A,F)$.

Recall that given an SH Lie pair $(L,A)$ with quotient space $B$, we get a Lie algebra object
$$
\mathfrak{B}=\Ha_{\mathrm{CE}}^\bullet(A, \mathbb{B})
$$
whose Lie bracket is induced by the Atiyah operator $\bm{\alpha}^B$ (Theorem~\ref{Cor:LieandLieModuleStr}).

Let $\Mod_{\mathfrak{B}}$ denote the category of $\mathfrak{B}$-modules. According to Theorem~\ref{Cor:LieandLieModuleStr} again, we are able to introduce
\begin{Def}
The Atiyah functor is 
$
\bm{A}: (E,\partial_A^E) \rightarrow \bigl(\Ha^\bullet_{\mathrm{CE}}(A,E), \bm{\alpha}^E\bigr)
$
from the category $\Mod_A$ of $A$-modules to the category $\Mod_{\mathfrak{B}}$ of $\mathfrak{B}$-modules. And for all $\phi \in \Hom_A(E,F)$, we have
\bd
 \bm{A}(\phi) = [\phi]: \Ha^\bullet_{\mathrm{CE}}(A,E) \rightarrow \Ha^\bullet_{\mathrm{CE}}(A,F).
\ed
\end{Def}

That $\bm{A}$ is well-defined relies on the following fact:~given any $\phi\in\Hom_A(E,F)$, the associated $[\phi]$ preserves the $\mathfrak{B}$-actions, i.e., for all $x\in \O(A)\otimes B$, $r\in \O(A)\otimes E$,
$$
\bm{\alpha}^F(x)(\phi(r))=\phi(\bm{\alpha}^E(x)r)+(\partial_A W)(x\otimes_{\O(A)} r),
$$
 for some $W\in \O(A)\otimes \Hom(B\otimes E,F)$.  In fact, we have $W=\delta(\phi)$. The proof of this fact is similar to that of Lemma~\ref{Lem:preAction}, and thus omitted.

\begin{Rm}
Inspired by Lemma~\ref{properties of Atiyah operators}, we may expect the Atiyah functor to enjoy  the following natural properties:
$$\bm{A}(E\otimes F)\cong \bm{A}(E)\otimes_{\mathcal{H}} \bm{A}(F),$$
$$\bm{A}(E^\vee)\cong \Hom_{\mathcal{H}}(\bm{A}(E), \mathcal{H}),$$
and
$$\bm{A}(\Hom(E,F))\cong \Hom_{\mathcal{H}}(\bm{A}(E), \bm{A}(F)).$$
However, some condition is needed to fulfill these isomorphisms. Further investigations of this question will be dealt with somewhere else.
\end{Rm}

\section{Invariance of Atiyah classes under infinitesimal deformations}\label{Invariance section}
In this section, let us fix an SH Lie pair $(L,A; Q_L \sim \lambda_\bullet)$ with the quotient space $B = L/A$, and an $A$-module $E$. We study infinitesimal deformations of the $L_\infty[1]$-structure $Q_L$ on $L$, and how the associated Atiyah classes $[\alpha^E]$ would be affected.

\subsection{Compatible infinitesimal deformations}
In what follows, $\hbar$ denotes a square zero formal parameter. An infinitesimal deformation, or a first order deformation, of the $L_{\infty}[1]$-algebra structure on $L$, namely that of $Q_L$, is a differential of the form
\bd
 Q(\hbar) = Q_L + \hbar Q_+:~~~\O(L)[\hbar] \rightarrow \O(L)[\hbar].
\ed
Here $Q_+$ is a degree $1$ derivation of $\O(L)$, and both $Q_L$ and $Q_+$ are $\k[\hbar]$-linear. It follows that
\bd
  [Q_L, Q_+] = Q_L\circ Q_+ + Q_+\circ Q_L=0.
\ed
In this circumstance,   $L[\hbar]$ has an $L_\infty[1]$-structure $Q(\hbar)$ which is deformed from $Q_L$.

As our motivation is to regard $L$ as a larger object extended from $A$, we only consider deformations of the following type:
\begin{Def}
  An infinitesimal deformation $Q(\hbar)$ of $Q_L$ is said to be $A$-compatible, if it is subject to the following two conditions:
  \begin{enumerate}
    \item[(1)] The $L_\infty[1]$-structure on $A$ is not deformed, i.e., the following diagram commutes:
    \bd
    \xymatrix{
     \O(L)[\hbar] \ar[r]^-{j^\vee} \ar[d]^-{Q(\hbar)}   & \O(A)[\hbar] \ar[d]^-{Q_A} \\
     \O(L)[\hbar] \ar[r]^-{j^\vee}   &    \O(A)[\hbar].
    }
    \ed
    \item[(2)] The $A$-module structure on $B$ is not deformed. This means the commutativity of
    \[  \xymatrix{
    A^\perp[\hbar]=B^\vee[\hbar] \ar[rr]^-{D^\perp} \ar@{^{(}->}[d] &    & \O(A)\otimes B^\vee[\hbar] \ar@{^{(}->}[d]\\
    L^\vee[\hbar]\ar[r]^-{Q(\hbar)}  & \O(L)[\hbar]\ar[r]^-{J} & \O(A) \otimes L^\vee[\hbar].
    }
    \]
   \end{enumerate}
\end{Def}

By choosing a splitting of Sequence~\eqref{SES de GVS}, so that $L\cong A\oplus B$ and that $\O(L)$ is identified with $\O(A)\otimes \O(B)$, the two compatible conditions are unraveled:~if $Q(\hbar) = Q_L + \hbar Q_+$, then the above
\begin{align*}
  \mbox{ Condition (1) } & \Leftrightarrow ~ Q_+( A^\vee)  \subset \O(A) \otimes \O^+(B),\\
  \mbox{ Condition (2) } & \Leftrightarrow ~ Q_+ (  B^\vee) \subset \O(A)  \otimes    \widehat{S}{^{\geq 2}(B^\vee)}.
\end{align*}

Similar to the decomposition of $Q_L$ in Equation~\eqref{T}, we denote the part of $Q_+$ that sends $A^\vee$ into $\O(A) \otimes B^\vee$ by $\delta_+$, the part that sends $A^\vee$ into $\O(A) \otimes  \widehat{S}{^{\geq 2}(B^\vee)}$ by $\mathcal{R}_+$, and the part that sends $B^\vee$ into $\O(A) \otimes S^i(B^\vee)$ by $T_+^i, i \geq 2$. Then the $A$-compatible infinitesimal deformation $Q(\hbar)$ has the form
\be\label{Eqt:Qhbarcomponents}
 Q(\hbar) = Q_L + \hbar\delta_+ + \hbar\mathcal{R}_+ + \hbar\sum_{i \geq 2}T_+^i.
\ee

\begin{Def}\label{Def:gaugeequivalence}
  \begin{enumerate}
    \item A gauge equivalence of $ \O(L)[\hbar] $ is an automorphism $\sigma = {1} + \hbar\lambda$ of the graded commutative algebra $\O(L)[\hbar]$, where $\lambda:\O(L)\rightarrow \O(L)$ is $\k$-linear, such that the following diagram commutes:
    \be\label{Digram:lambdacommutes}
    \xymatrix{
    \O(B)[\hbar] \ar[r] \ar[d]^-{1} & \O(L)[\hbar] \ar[r]^-{j^\vee} \ar[d]^-{\sigma=1+\hbar\lambda}   & \O(A)[\hbar] \ar[d]^-{1} \\
    \O(B)[\hbar] \ar[r] & \O(L)[\hbar] \ar[r]^-{j^\vee}   & \O(A)[\hbar].
    }
    \ee
    \item Two $A$-compatible infinitesimal deformations $Q(\hbar)$ and $\bar{Q}(\hbar)$ of   $Q_L$ are said to be gauge equivalent if there exists a gauge equivalence $\sigma=1+\hbar\lambda  $ such that the following diagram commutes:
    \bd
    \xymatrix{
    \O(L)[\hbar] \ar[rrr]^-{\sigma=1+\hbar\lambda} \ar[d]^-{Q(\hbar)}  && & \O(L)[\hbar] \ar[d]^-{\bar{Q}(\hbar)} \\
    \O(L)[\hbar] \ar[rrr]^-{\sigma=1+\hbar\lambda} &&  & \O(L)[\hbar],
    }
  \ed
  i.e., $\sigma$ is an isomorphism of $L_\infty[1]$-algebras $(L[\hbar],Q[\hbar])\cong (L[\hbar],\bar{Q}[\hbar])$.
  \end{enumerate}
\end{Def}
Assume that $Q(\hbar)$ and $\bar{Q}(\hbar)$ are connected by the gauge equivalence $\sigma = 1+\hbar\lambda$. Since $\sigma$ is an algebra automorphism, it follows that $\lambda$ is a degree $0$ derivation of $\O(L)$. Note that $\sigma^{-1}=1-\hbar\lambda$. It follows from a simple computation that
$$
Q(\hbar)-\bar{Q}(\hbar)=\hbar[Q_L,\lambda].
$$

Recall that $\O(L)=\O(A)\otimes \O(B)$. The commutative property of Diagram~\eqref{Digram:lambdacommutes} implies that we can write
\be\label{Eqt:lambdadecompose}
 \lambda = \sum_{k\geq1} \Psi_k,\;\;\text{where}\;\;\; \Psi_k:~A^\vee \rightarrow \O(A) \otimes S^k(B^\vee).
\ee
All these $\Psi_k$ are treated as degree $0$ derivations of $\O(L)$ which act trivially on $B^\vee$.

\subsection{Gauge invariance of Atiyah classes}
Let $Q(\hbar) = Q_L + \hbar Q_+$ be an $A$-compatible infinitesimal deformation of $Q_L$. Consider the associated Atiyah cocycle $\alpha^{E[\hbar]}$ of the SH Lie pair $(L[\hbar],A[\hbar])$ with respect to the $A[\hbar]$-module $E[\hbar]$. By Lemma~\ref{Important},
\bd
\alpha^{E[\hbar]} = [\partial_{A}, \delta + \hbar\delta_+] = \alpha^E + \hbar[\partial_A ,\delta_+],
\ed
where $\delta$ and $\delta_+$ are the components of $Q_L$ and $Q_+$ specified respectively in Equations \eqref{T} and \eqref{Eqt:Qhbarcomponents}.

The main result in this section is the gauge invariance of the Atiyah class $[\alpha^{E[\hbar]}]$.
\begin{Thm}\label{Thm:gaugeinvariance}
Let $Q(\hbar) = Q_L + \hbar Q_+$ and $\bar{Q}(\hbar) = Q_L + \hbar \bar{Q}_+$ be two gauge equivalent $A$-compatible infinitesimal deformations of $Q_L$. Then the associated Atiyah classes coincide:
    \bd
  [\alpha^{E[\hbar]}] = [\overline{\alpha^{E[\hbar]}}] \in H^2(A[\hbar],(B[\hbar])^\vee \otimes \End(E[\hbar])) \cong H^2(A,B^\vee\otimes\End(E))[\hbar].
  \ed
\end{Thm}
\bp
Let the gauge equivalence $\sigma$ be as in Definition \ref{Def:gaugeequivalence}. The $\lambda$ operator is defined in Equation~\eqref{Eqt:lambdadecompose}. Further assume that
\begin{align*}
 Q(\hbar) &= Q_A + \delta + D^\perp + \mathcal{R} + \sum_{j \geq 2}T^j + \hbar\delta_+ + \hbar\mathcal{R}_+ + \hbar\sum_{i \geq 2}T_+^i,\\
 \bar{Q}(\hbar) &= Q_A + \delta + D^\perp + \mathcal{R} + \sum_{j \geq 2}T^j + \hbar\bar{\delta}_+ + \hbar\bar{\mathcal{R}}_+ + \hbar\sum_{i \geq 2}\bar{T}_+^i,
\end{align*}
are explained as earlier. Applying the equation
\bd
 \sigma \circ Q(\hbar) = \bar{Q}(\hbar) \circ \sigma:~~\O(L)[\hbar] \rightarrow \O(L)[\hbar]
\ed
to an element $\xi \in A^\vee$, we have
\begin{align*}
  \sigma &(Q(\hbar)(\xi)) = \left(1 + \hbar\sum_k\Psi_k\right)(Q_A(\xi) + \delta(\xi) + \mathcal{R}(\xi) + \hbar\delta_+(\xi) + \hbar\mathcal{R}_+(\xi)) \\
  &= Q_A(\xi) + \delta(\xi) + \mathcal{R}(\xi) + \hbar\delta_+(\xi) + \hbar\mathcal{R}_+(\xi) + \hbar\sum_k\Psi_k(Q_A(\xi) + \delta(\xi) + \mathcal{R}(\xi)),
\end{align*}
and
\begin{align*}
  &\bar{Q}(\hbar)(\sigma(\xi)) = \left(Q_A + \delta + D^\perp + \mathcal{R} + \sum_{j \geq 2}T^j + \hbar\bar{\delta}_+ + \hbar\bar{\mathcal{R}}_+ + \hbar\sum_{i \geq 2}\bar{T}_+^i\right)\left(\xi + \hbar\sum_k\Psi_k(\xi)\right) \\
  &= Q_A(\xi) + \delta(\xi) + \mathcal{R}(\xi) + \hbar\bar{\delta}_+(\xi) + \hbar\bar{\mathcal{R}}_+(\xi) + \hbar\sum_k \left(Q_A + \delta + D^\perp + \mathcal{R} + \sum_{j \geq 2}T^j\right)(\Psi_k(\xi)).
\end{align*}
Comparing the $\hbar \O(A)\otimes B^\vee$-component of both sides, one gets
\bd
 \Psi_1(Q_A(\xi)) + \delta_+(\xi) = \bar{\delta}_+(\xi) + Q_A(\Psi_1(\xi)) + D^\perp(\Psi_1(\xi)),
\ed
which implies that
\bd
 \delta_+ - \bar{\delta}_+ = Q_A \circ \Psi_1 - \Psi_1 \circ Q_A + D^\perp \circ \Psi_1:~ A^\vee \rightarrow \O(A) \otimes B^\vee.
\ed

Hence, we have, for all $e \in E$,
\begin{align*}
  &\quad\alpha^{E[\hbar]}(e) - \overline{\alpha^{E[\hbar]}}(e) = \hbar[\partial_A^E, \delta_+ - \bar{\delta}_+](e) = \hbar(\delta_+ - \bar{\delta}_+) (\partial_A^E(e)) \\
  &= \hbar(Q_A \circ \Psi_1 - \Psi_1 \circ Q_A + D^\perp \circ \Psi_1)\left(D^E(e)\right) \\
  &= \hbar(Q_A \circ \Psi_1 \circ D^E + \Psi_1 \circ D^E \circ D^E + D^\perp \circ \Psi_1 \circ D^E)(e).\;\quad\quad(\text{by Equation~\eqref{MC-EQ}})
\end{align*}
Here $D^E:E\rightarrow \O(A)\otimes E$ defines the $A$-module structure on $E$.

Now let
\bd
W = [\Psi_1,D^E]:~\O(A) \otimes E \rightarrow \O(A) \otimes B^\vee \otimes E
\ed
be the graded commutator of $\Psi_1$ and $D^E$. One easily finds it is an $\O(A)$-linear map.

Then we have
\begin{align*}
  \alpha^{E[\hbar]}(e) - \overline{\alpha^{E[\hbar]}}(e) &= \hbar((Q_A + D^\perp + D^E)((\Psi_1 \circ D^E)(e)) + (W \circ D^E)(e)) \\
  &= \hbar(\partial_A  \circ W + W \circ \partial_A)(e) = \hbar \partial_A W(e),
\end{align*}
which implies that
\bd
 \alpha^{E[\hbar]} - \overline{\alpha^{E[\hbar]}} = \hbar\partial_A W.
\ed
  This completes the proof.
\ep

\section{Appendix: Morphisms of SH Lie algebras}\label{appendix}
We prove the equivalence of the two definitions of morphisms of SH Lie algebras as in Definition~\ref{morphism of SH Lie algebras}.

Let $\phi: \O(L^\prime) \rightarrow \O(L)$ be a morphism of $\k$-algebras such that
\bd
 \phi \circ Q_{L^\prime} = Q_L \circ \phi: \O(L^\prime) \rightarrow \O(L).
\ed
Assume that $\phi = \sum_{k\geq0} \phi_k$, where $\phi_k: (L^\prime)^\vee \rightarrow S^k(L^\vee)$.  Define a family of degree zero linear maps
\bd
  f_k=(-1)^{k+1} \phi_k^\vee:~S^k(L) \rightarrow L^\prime, \;\; k \geq 0.
\ed
We show that $\{f_k\}$ satisfies the two requirements as in Definition~\ref{morphism of SH Lie algebras}.

By assumption, we have
\be\label{morphism equation}
\phi \circ Q_{L^\prime}(\xi) = Q_L \circ \phi(\xi) \in \O(L)
\ee
for all homogeneous $\xi \in (L^\prime)^\vee$. Note that
\begin{align*}
  &\quad \text{LHS of Equation~\eqref{morphism equation}}  = \phi\left(\langle \lambda_0^\prime,\xi \rangle + Q_1^\prime(\xi) + \cdots \right) = \langle \lambda_0^\prime,\xi \rangle + \sum_{n \geq 1}\langle (\phi_0)^{\odot n}, Q_n^\prime(\xi) \rangle + \cdots;
\end{align*}
and
\begin{align*}
  &\quad \text{RHS of Equation~\eqref{morphism equation}} = Q_L\left(\langle \phi_0, \xi \rangle + \phi_1(\xi) + \cdots \right) = Q_0(\phi_1(\xi)) + \cdots = \langle \lambda_0, \phi_1(\xi) \rangle + \cdots.
\end{align*}
Comparing the $\k = S^0(L^\vee)$-component of both sides, one gets
\begin{align*}
  \langle \lambda_0^\prime,\xi \rangle + \sum_{n \geq 1}\langle \phi_0^{\odot n}, Q_n^\prime(\xi) \rangle = \langle \lambda_0, \phi_1(\xi) \rangle,
\end{align*}
which is equivalent to
\bd
\langle \lambda_0^\prime + \sum_{n \geq 1}\frac{1}{n!}(-1)^n\lambda_n^\prime(\phi_0,\cdots,\phi_0), \xi \rangle = \langle \lambda_0, \phi_1(\xi) \rangle,
\ed
or Equation~\eqref{f_0-relation}.

We further investigate the $S^n(L^\vee)(n \geq 1)$-component of Equation~\eqref{morphism equation}. For all $u_i \in L, i=1,\cdots,n$, we have
\begin{align*}
   &\quad \text{LHS} = \langle \phi \circ Q_{L^\prime}(\xi), u_1 \odot \cdots \odot u_n \rangle \\
   &= \langle Q_{L^\prime}(\xi), \sum_{\substack{i_1,\cdots,i_r \geq 1 \\ i_1 + \cdots + i_r = n}}\sum_{\tau \in \sh(i_1,\cdots,i_r)}\sum_{j \geq 0}\epsilon(\tau) (-1)^{n+r+j} \frac{1}{r!} f_0^{\odot j}\odot f_{i_1}\odot\cdots \odot f_{i_r}(u_{\tau(1)},\cdots,u_{\tau(n)}) \rangle \\
   &= \langle \xi, \sum_{\substack{i_1,\cdots,i_r \geq 1 \\ i_1 + \cdots + i_r = n}}\sum_{\tau \in \sh(i_1,\cdots,i_r)}\sum_{j \geq 0}\epsilon(\tau) (-1)^{\abs{\xi}+n} \frac{1}{(r+j)!} \lambda_{r+j}^\prime(f_0,\cdots,f_0,f_{i_1}(\cdots),\cdots,f_{i_r}(\cdots))\rangle
\end{align*}
and
\begin{align*}
  \quad \text{RHS} &= \langle Q_L \circ \phi(\xi), u_1 \odot \cdots \odot u_n \rangle \\
  &= \langle \phi(\xi), \sum_{\substack{k,l \geq 0 \\ k+l =n}}\sum_{\sigma \in \sh(l,k)}\epsilon(\sigma)(-1)^{\abs{\xi}+l} \lambda_l(u_{\sigma(1)},\cdots, u_{\sigma(l)})\odot\cdots\odot u_{\sigma(n)} \rangle \\
  &= \langle \xi, \sum_{\substack{k,l \geq 0 \\ k+l =n}}\sum_{\sigma \in \sh(l,k)}\epsilon(\sigma)(-1)^{\abs{\xi}+n} f_{k+1}( \lambda_l(u_{\sigma(1)}, \cdots, u_{\sigma(l)}),\cdots,u_{\sigma(n)}) \rangle.
\end{align*}
Thus Equation~\eqref{morphism relation} also holds once we assume Equation~\eqref{Q-morphism}.

The inverse implication ``Equations~\eqref{f_0-relation}+\eqref{morphism relation} $\Longrightarrow$ Equation~\eqref{morphism equation}'' is also clear from the previous argument. This completes the proof.

%\begin{Rm}
%The ordinary category of Lie algebras is not a full sub-category of $L_\infty[1]$-algebras. A morphism of Lie algebras in the usual sense is a $\k$-linear map of the underlying graded vector spaces preserving Lie brackets, while a morphism of $L_\infty[1]$-algebras is a smooth map of the underlying graded vector spaces preserving  homological vector fields.
%\end{Rm}

\subsection*{Acknowledgments}
We would like to thank Camille Laurent-Gengoux, Xiaobo Liu and Ping Xu for useful discussions and comments. The authors are indebted to an anonymous referee who provided many suggestions which helped us a lot to improve an early version of this paper.
Zhuo Chen is also grateful to Sheffield University for its hospitality.
Honglei Lang and Maosong Xiang also would like to thank the Department of Mathematics at Penn State for its hospitality and China Scholarship Council for the funding during their visit at Penn State.


\begin{bibdiv}
  \begin{biblist}
	
  \bib{AC}{article}{
   author={Abad, Camilo Arias},
   author={Crainic, Marius},
   title={Representations up to homotopy of Lie algebroids},
   journal={J. Reine Angew. Math.},
   volume={663},
   date={2012},
   pages={91--126},
}

  \bib{AKSZ}{article}{
   author={Alexandrov, M.},
   author={Kontsevich, M.},
   author={Schwarz, A.},
   author={Zaboronsky, O.},
   title={The geometry of the master equation and topological quantum field
   theory},
   journal={Internat. J. Modern Phys. A},
   volume={12},
   date={1997},
   number={7},
   pages={1405--1429},
}

  \bib{Atiyah}{article}{
   author={Atiyah, M. F.},
   title={Complex analytic connections in fibre bundles},
   journal={Trans. Amer. Math. Soc.},
   volume={85},
   date={1957},
   pages={181--207},
}


  % \bib{Bru}{article}{
  % author={Bruce, Andrew James},
  % title={From $L_\infty$-algebroids to higher Schouten/Poisson structures},
  % eprint={arXiv:~1007.1389v3}
%}

  \bib{BC}{article}{
   author={Baez, John C.},
   author={Crans, Alissa S.},
   title={Higher-dimensional algebra. VI. Lie 2-algebras},
   journal={Theory Appl. Categ.},
   volume={12},
   date={2004},
   pages={492--538},
}

  \bib{BV}{article}{
   author={Bashkirov, Denis},
   author={Voronov, Alexander A.},
   title={The BV formalism for $L_\infty$-algebras},
   journal={J. Homotopy Relat. Struct.},
   volume={12},
   date={2017},
   number={2},
   pages={305--327},
   issn={2193-8407},
%   review={\MR{3654355}},
%   doi={10.1007/s40062-016-0129-z},
}

   \bib{Bordemann}{article}{
   author={Bordemann, Martin},
   title={Atiyah classes and equivariant connections on homogeneous spaces},
   conference={
      title={Travaux math\'ematiques. Volume XX},
   },
   book={
      series={Trav. Math.},
      volume={20},
      publisher={Fac. Sci. Technol. Commun. Univ. Luxemb., Luxembourg},
   },
   date={2012},
   pages={29--82},
   %review={\MR{3014184}},
}

   \bib{Bottacin}{article}{
   author={Bottacin, Francesco},
   title={Atiyah classes of Lie algebroids},
   conference={
      title={Current trends in analysis and its applications},
   },
   book={
      series={Trends Math.},
      publisher={Birkh\"auser/Springer, Cham},
   },
   date={2015},
   pages={375--393},
%   review={\MR{3496770}},
}

  \bib{Calaque}{article}{
   author={Calaque, Damien},
   title={A PBW theorem for inclusions of (sheaves of) Lie algebroids},
   journal={Rend. Semin. Mat. Univ. Padova},
   volume={131},
   date={2014},
   pages={23--47},
   %issn={0041-8994},
   %review={\MR{3217749}},
   %doi={10.4171/RSMUP/131-3},
}


   \bib{CCT}{article}{
   author={Calaque, Damien},
   author={C{\u{a}}ld{\u{a}}raru, Andrei},
   author={Tu, Junwu},
   title={PBW for an inclusion of Lie algebras},
   journal={J. Algebra},
   volume={378},
   date={2013},
   pages={64--79},
  % issn={0021-8693},
  % review={\MR{3017014}},
  % doi={10.1016/j.jalgebra.2012.12.008},
}

 \bib{CV}{article}{
   author={Calaque, Damien},
   author={Van den Bergh, Michel},
   title={Hochschild cohomology and Atiyah classes},
   journal={Adv. Math.},
   volume={224},
   date={2010},
   number={5},
   pages={1839--1889},
  % issn={0001-8708},
  % review={\MR{2646112}},
  % doi={10.1016/j.aim.2010.01.012},
}

   \bib{CSX14}{article}{
   author={Chen, Zhuo},
   author={Sti{\'e}non, Mathieu},
   author={Xu, Ping},
   title={A Hopf algebra associated with a Lie pair},
   language={English, with English and French summaries},
   journal={C. R. Math. Acad. Sci. Paris},
   volume={352},
   date={2014},
   number={11},
   pages={929--933},
   %issn={1631-073X},
   %review={\MR{3268765}},
   %doi={10.1016/j.crma.2014.09.010},
}

   \bib{CSX}{article}{
   author={Chen, Zhuo},
   author={Sti{\'e}non, Mathieu},
   author={Xu, Ping},
   title={From Atiyah classes to homotopy Leibniz algebras},
   journal={Comm. Math. Phys.},
   volume={341},
   date={2016},
   number={1},
   pages={309--349},
}

   \bib{Costello}{article}{
   author={Costello, Kevin},
   title={A geometric construction of the Witten genus, I},
   conference={
      title={Proceedings of the International Congress of Mathematicians.
      Volume II},
   },
   book={
      publisher={Hindustan Book Agency, New Delhi},
   },
   date={2010},
   pages={942--959},
  % review={\MR{2827826}},
}


%   \bib{GM}{book}{
%   author={Gelfand, Sergei I.},
%   author={Manin, Yuri I.},
%   title={Methods of homological algebra},
%   series={Springer Monographs in Mathematics},
%   edition={2},
%   publisher={Springer-Verlag, Berlin},
%   date={2003},
%}

  \bib{Grinberg}{article}{
   author={Grinberg, Darij},
   title={Poincar\'{e}-Birkhoff-Witt type results for inclusions of Lie algebras},
   note={Master's thesis--Massachusetts Institute of Technology, Boston, MA},
   date={2011},
}

   \bib{Kap}{article}{
   author={Kapranov, M.},
   title={Rozansky-Witten invariants via Atiyah classes},
   journal={Compositio Math.},
   volume={115},
   date={1999},
   number={1},
   pages={71--113},
}

   \bib{Kon1}{article}{
   author={Kontsevich, Maxim},
   title={Operads and motives in deformation quantization},
   note={Mosh\'e Flato (1937--1998)},
   journal={Lett. Math. Phys.},
   volume={48},
   date={1999},
   number={1},
   pages={35--72},
   %issn={0377-9017},
   %review={\MR{1718044 (2000j:53119)}},
   %doi={10.1023/A:1007555725247},
}

   \bib{Kon2}{article}{
   author={Kontsevich, Maxim},
   title={Deformation quantization of Poisson manifolds},
   journal={Lett. Math. Phys.},
   volume={66},
   date={2003},
   number={3},
   pages={157--216},
   %issn={0377-9017},
   %review={\MR{2062626 (2005i:53122)}},
   %doi={10.1023/B:MATH.0000027508.00421.bf},
}



   \bib{KS}{article}{
   author={Kontsevich, Maxim},
   author={Soibelman, Yan},
   title={Deformation theory. I. Draft of a book},
   eprint={http://www.math.ksu.edu/~soibel/},
   date={2006}
}


   \bib{LM}{article}{
   author={Lada, Tom},
   author={Markl, Martin},
   title={Strongly homotopy Lie algebras},
   journal={Comm. Algebra},
   volume={23},
   date={1995},
   number={6},
   pages={2147--2161},
}


   \bib{LS}{article}{
   author={Lada, Tom},
   author={Stasheff, Jim},
   title={Introduction to SH Lie algebras for physicists},
   journal={Internat. J. Theoret. Phys.},
   volume={32},
   date={1993},
   number={7},
   pages={1087--1103},
}

   \bib{LSX1}{article}{
   author={Laurent-Gengoux, Camille},
   author={Sti{\'e}non, Mathieu},
   author={Xu, Ping},
   title={Exponential map and $L_\infty$ algebra associated to a Lie
   pair},
   language={English, with English and French summaries},
   journal={C. R. Math. Acad. Sci. Paris},
   volume={350},
   date={2012},
   number={17-18},
   pages={817--821},
}

   \bib{LSX2}{article}{
   author={Laurent-Gengoux, Camille},
   author={Sti{\'e}non, Mathieu},
   author={Xu, Ping},
   title={Kapranov dg-manifolds and Poincar\'{e}-Birkhoff-Witt isomorphisms},
   eprint={arXiv:~1408.2903},
}

   \bib{LV}{article}{
   author={Laurent-Gengoux, Camille},
   author={Voglaire, Y.},
   title={Invariant connections and PBW theorem for Lie groupoid pairs},
   eprint={arXiv:~1507.01051},
   date={2015},
}


   \bib{Mehta}{article}{
   author={Mehta, Rajan Amit},
   title={Lie algebroid modules and representations up to homotopy},
   journal={Indag. Math. (N.S.)},
   volume={25},
   date={2014},
   number={5},
   pages={1122--1134},
   %issn={0019-3577},
   %review={\MR{3264789}},
   %doi={10.1016/j.indag.2014.07.013},
}

   \bib{MSX}{article}{
   author={Mehta, Rajan Amit},
   author={Sti{\'e}non, Mathieu},
   author={Xu, Ping},
   title={The Atiyah class of a dg-vector bundle},
   language={English, with English and French summaries},
   journal={C. R. Math. Acad. Sci. Paris},
   volume={353},
   date={2015},
   number={4},
   pages={357--362},
  % issn={1631-073X},
  % review={\MR{3319134}},
  % doi={10.1016/j.crma.2015.01.019},
}

  \bib{Molino1}{article}{
   author={Molino, Pierre},
   title={Classe d'Atiyah d'un feuilletage et connexions transverses
   projetables. },
   language={French},
   journal={C. R. Acad. Sci. Paris S\'er. A-B},
   volume={272},
   date={1971},
   pages={A779--A781},
}

   \bib{Molino2}{article}{
   author={Molino, Pierre},
   title={Propri\'et\'es cohomologiques et propri\'et\'es topologiques des
   feuilletages \`a connexion transverse projetable},
   language={French},
   journal={Topology},
   volume={12},
   date={1973},
   pages={317--325},
}

\bib{Nguyen}{article}{
   author={Nguyen-van, Hai},
   title={Relations entre les diverses obstructions relatives \`a
   l'existence d'une connexion lin\'eaire invariante sur un expace
   homog\`ene},
   language={French},
   journal={C. R. Acad. Sci. Paris},
   volume={260},
   date={1965},
   pages={45--48},
  %review={\MR{0176417}},
}

   \bib{Quillen}{article}{
   author={Quillen, Daniel},
   title={Superconnections and the Chern character},
   journal={Topology},
   volume={24},
   date={1985},
   number={1},
   pages={89--95},
}


   \bib{RW}{article}{
   author={Rozansky, L.},
   author={Witten, E.},
   title={Hyper-K\"ahler geometry and invariants of three-manifolds},
   journal={Selecta Math. (N.S.)},
   volume={3},
   date={1997},
   number={3},
   pages={401--458},
  % issn={1022-1824},
  %review={\MR{1481135 (98m:57041)}},
  %doi={10.1007/s000290050016},
}

   \bib{Stasheff}{article}{
   author={Stasheff, Jim},
   title={Higher homotopy structures:~then and now},
   note={The talk at Opening workshop of Higher Structures in Geometry and Physics, MPI Bonn 2016},
   %eprint={https://ncatlab.org/nlab/files/StasheffHomotopyStructuresReview.pdf},
   }

  \bib{Shoikhet}{article}{
   author={Shoikhet, Boris},
   title={On the Duflo formula for $L_\infty$-algebras and $Q$-manifolds},
   eprint={https://arxiv.org/abs/math/9812009},
   date={1998},
}

   \bib{Luca}{article}{
   author={Vitagliano, Luca},
   title={Representations of homotopy Lie-Rinehart algebras},
   journal={Math. Proc. Cambridge Philos. Soc.},
   volume={158},
   date={2015},
   number={1},
   pages={155--191},
   %issn={0305-0041},
   %review={\MR{3300319}},
   %doi={10.1017/S0305004114000541},
}

   \bib{Vor}{article}{
   author={Voronov, Theodore},
   title={Higher derived brackets and homotopy algebras},
   journal={J. Pure Appl. Algebra},
   volume={202},
   date={2005},
   number={1-3},
   pages={133--153},
}

   \bib{Vor2}{article}{
   author={Voronov, Theodore},
   title={Higher derived brackets for arbitrary derivations},
   conference={
      title={Travaux math\'ematiques. Fasc. XVI},
   },
   book={
      series={Trav. Math., XVI},
      publisher={Univ. Luxemb., Luxembourg},
   },
   date={2005},
   pages={163--186},
}

   \bib{Wang}{article}{
   author={Wang, Hsien-chung},
   title={On invariant connections over a principal fibre bundle},
   journal={Nagoya Math. J.},
   volume={13},
   date={1958},
   pages={1--19},
   %issn={0027-7630},
   %review={\MR{0107276}},
}

  \end{biblist}
\end{bibdiv}
\end{document}